\documentclass[a4paper,12pt]{article}
\usepackage{graphicx}
\usepackage{psfrag}

\usepackage{mathrsfs}
\usepackage{amssymb}
\usepackage{amsfonts}
\usepackage{latexsym}
\usepackage{amsmath}
\usepackage{amsthm}
\usepackage[cp1251]{inputenc}
\usepackage[english]{babel}

\newcommand{\er}[1]{\textrm{(\ref{#1})}}
\def\lb{\label}

\theoremstyle{plain}

\newtheorem{theorem}{\bf Theorem}[section]
\newtheorem{lemma}[theorem]{\bf Lemma}

\theoremstyle{remark}

\renewcommand{\a}{\alpha}            	  
                 
\newcommand{\g}{\gamma}                  
\newcommand{\G}{\Gamma}

\newcommand{\D}{\Delta}                
\newcommand{\ve}{\varepsilon}          
            \newcommand{\cH}{\mathcal{H}}     
\newcommand{\e}{\eta}                  
\newcommand{\vt}{\vartheta}            
\newcommand{\vT}{\Theta}               
                
\renewcommand{\l}{\lambda}          \newcommand{\cM}{\mathcal{M}}     
               
\newcommand{\m}{\mu}                   
\newcommand{\n}{\nu}              \newcommand{\cP}{\mathcal{P}}     
\renewcommand{\r}{\rho}                  
\newcommand{\s}{\sigma}           \newcommand{\cR}{\mathcal{R}}      

           \newcommand{\cS}{\mathcal{S}}     
\renewcommand{\t}{\tau}             \newcommand{\cT}{\mathcal{T}}     
\newcommand{\f}{\phi}                  
\newcommand{\F}{\Phi}                  
\newcommand{\vp}{\varphi}              
                   
\newcommand{\p}{\psi}                   
             \newcommand{\cZ}{\mathcal{Z}}      		
\renewcommand{\o}{\omega}

\newcommand{\x}{\xi}

  \def\mA{{\mathscr A}}
 \def\mB{{\mathscr B}}

  \def\mH{{\mathscr H}}

  \def\mU{{\mathscr U}}

\newcommand{\gD}{\mathfrak{D}}

\newcommand{\gH}{\mathfrak{H}}

\newcommand{\gR}{\mathfrak{R}}
\newcommand{\gS}{\mathfrak{S}}

\def\Z{\mathbb{Z}}
\def\R{\mathbb{R}}
\def\C{\mathbb{C}}

\def\N{\mathbb{N}}

\def\qqq{\qquad}
\def\qq{\quad}
\let\ge\geqslant
\let\le\leqslant
\let\geq\geqslant
\let\leq\leqslant
\newcommand{\ca}{\begin{cases}}
\newcommand{\ac}{\end{cases}}
\newcommand{\ma}{\begin{pmatrix}}
\newcommand{\am}{\end{pmatrix}}

\def\lt{\biggl}
\def\rt{\biggr}
\def\lra{\Leftrightarrow}
\renewcommand{\[}{\begin{equation}}
\renewcommand{\]}{\end{equation}}
\def\wt{\widetilde}
\def\pa{\partial}
\def\sm{\setminus}
\def\es{\emptyset}
\def\no{\noindent}
\def\ol{\overline}
\def\iy{\infty}
\def\ev{\equiv}
\def\/{\over}
\def\ts{\times}

\def\os{\oplus}
\def\ss{\subset}

\def\wh{\widehat}

\def\Im{\mathop{\rm Im}\nolimits}
\def\supp{\mathop{\rm supp}\nolimits}

\def\Tr{\mathop{\rm Tr}\nolimits}
\def\const{\mathop{\rm const}\nolimits}

\def\BBox{\hspace{1mm}\vrule height6pt width5.5pt depth0pt \hspace{6pt}}

\begin{document}
\title{Zigzag periodic nanotube in
magnetic field}
\author{
Evgeny Korotyaev
\begin{footnote} {
Institut f\"ur  Mathematik,  Humboldt Universit\"at zu Berlin,
Rudower Chaussee 25, 12489 Berlin, Germany,
e-mail: evgeny@math.hu-berlin.de\ \ 
}
\end{footnote}
\and Igor Lobanov
\begin{footnote} {
Mathematical Faculty, Mordovian State University,
430000 Saransk, Russia,
e-mail: lobanov@math.mrsu.ru }
\end{footnote}
}

\maketitle

\begin{abstract}
\no We consider the magnetic Schr\"odinger operator on the so-called zigzag periodic metric graph (a quasi 1D continuous model of zigzag nanotubes) with a periodic potential. The magnetic field (with the amplitude $B\in \R$) is uniform and it is parallel to the axis of the  nanotube.
The spectrum of this operator consists of an absolutely continuous part (spectral bands separated by gaps) plus an infinite number of eigenvalues  with infinite multiplicity. We describe all compactly supported eigenfunctions with the same eigenvalue. We define a Lyapunov function,
which is analytic on some Riemann surface. On each
sheet, the Lyapunov function has the same properties 
as in the scalar case, but it has branch points, which we call resonances. We prove that all  resonances are real. We determine the asymptotics of the periodic and anti-periodic spectrum and of the resonances at high energy. 
We show that endpoints of the gaps are periodic or 
anti-periodic eigenvalues or resonances (real branch points of the Lyapunov function). We describe the spectrum as functions of $B$.
For example, if $B\to B_{k,m}={4({\pi\/2}-{\pi k\/N}+\pi m)\/\sqrt 3 \cos {\pi\/2N}}, k=1,2,..,N, m\in \Z$, then some spectral band shrinkes into a flat band i.e., an eigenvalue of infinite multiplicity.

\end{abstract}

\section{Introduction  and main results}
\setcounter{equation}{0}

Consider the Schr\"odinger operator $\mH =(-i\nabla-\mA)^2+q$ on the so-called zigzag periodic metric graph $\G^{(N)},N\ge 1$ 
with a periodic potential $q$ and with a uniform magnetic field
$\mB=B(0,0,1)\in \R^3,\ B\in\R$. The corresponding vector potential 
$\mA$ is given by
\[
\lb{vep}
\qqq
\mA(x)={1\/2}[\mB,x]={B\/2}\left(-x_2,x_1,0\right),\quad x=(x_1,x_2,x_3)\in\R^3.
\]
%
\begin{figure}\lb{fig1}
\centering
\noindent
(a){
\tiny
\psfrag{g001}[l][l]{$\Gamma_{0,0,1}$}
\psfrag{g002}[l][l]{$\Gamma_{0,0,2}$}
\psfrag{g003}[l][l]{$\Gamma_{0,0,N}$}
\psfrag{g011}[c][c]{$\Gamma_{0,1,1}$}
\psfrag{g012}[c][c]{$\Gamma_{0,1,2}$}
\psfrag{g013}[c][c]{$\Gamma_{0,1,N}$}
\psfrag{g021}[c][c]{$\Gamma_{0,2,1}$}
\psfrag{g022}[c][c]{$\Gamma_{0,2,2}$}
\psfrag{g023}[c][c]{$\Gamma_{0,2,N}$}
\psfrag{g-101}[c][c]{$\Gamma_{-1,0,1}$}
\psfrag{g-102}[l][l]{$\Gamma_{-1,0,2}$}
\psfrag{g-103}[l][l]{$\Gamma_{-1,0,N}$}
\psfrag{g-111}[c][c]{$\Gamma_{-1,1,1}$}
\psfrag{g-112}[c][c]{$\Gamma_{-1,1,2}$}
\psfrag{g-113}[c][c]{$\Gamma_{-1,1,N}$}
\psfrag{g-121}[c][c]{$\Gamma_{-1,2,1}$}
\psfrag{g-122}[c][c]{$\Gamma_{-1,2,2}$}
\psfrag{g-123}[c][c]{$\Gamma_{-1,2,N}$}
\psfrag{g101}[l][l]{$\Gamma_{1,0,1}$}
\psfrag{g102}[l][l]{$\Gamma_{1,0,2}$}
\psfrag{g103}[l][l]{$\Gamma_{1,0,N}$}
\psfrag{g111}[c][c]{$\Gamma_{1,1,1}$}
\psfrag{g112}[c][c]{$\Gamma_{1,1,2}$}
\psfrag{g113}[c][c]{$\Gamma_{1,1,N}$}
\psfrag{g121}[c][c]{$\Gamma_{1,2,1}$}
\psfrag{g122}[c][c]{$\Gamma_{1,2,2}$}
\psfrag{g123}[c][c]{$\Gamma_{1,2,N}$}
\includegraphics[width=.65\textwidth,height=.5\textwidth]{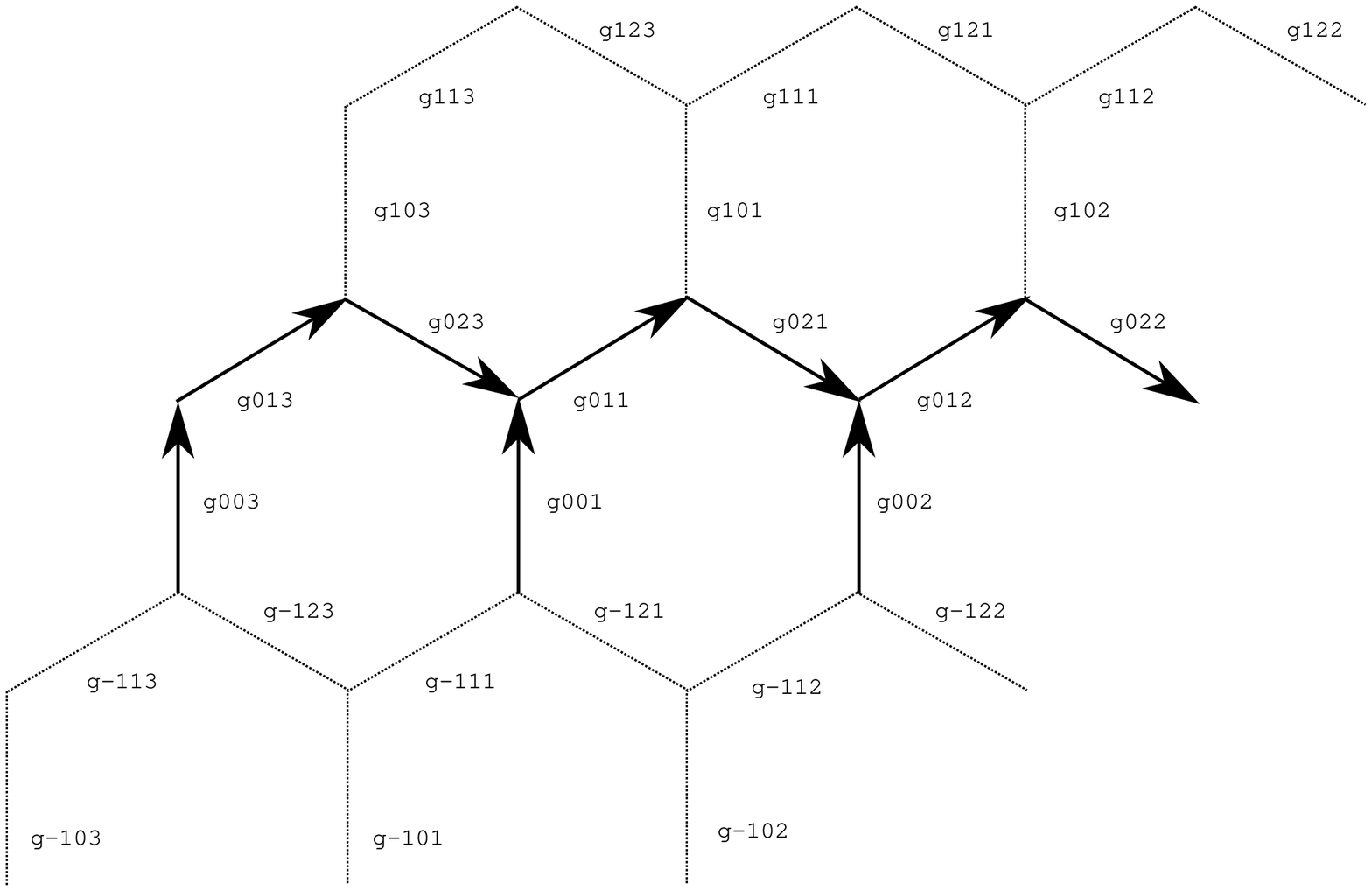}
}
(b){
\tiny
\psfrag{g-10}[l][l]{$\Gamma_{-1,0}$}
\psfrag{g-11}[l][l]{$\Gamma_{-1,1}$}
\psfrag{g-12}[l][l]{$\Gamma_{-1,2}$}
\psfrag{g00}[l][l]{$\Gamma_{0,0}$}
\psfrag{g01}[l][l]{$\Gamma_{0,1}$}
\psfrag{g02}[l][l]{$\Gamma_{0,2}$}
\psfrag{g10}[c][c]{$\Gamma_{1,0}$}
\psfrag{g11}[c][c]{$\Gamma_{1,2}$}
\psfrag{g12}[c][c]{$\Gamma_{1,3}$}
\includegraphics[height=.5\textwidth]{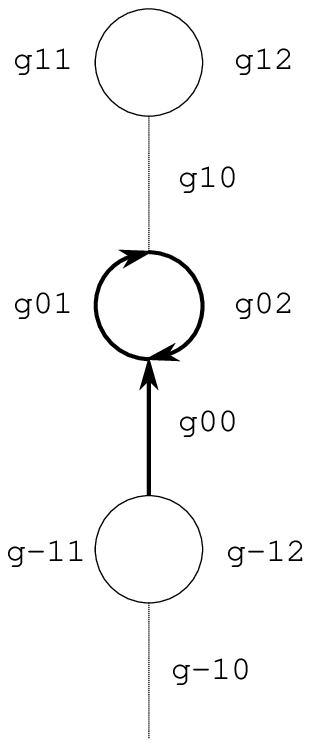}
}
\caption{Zigzag graph $\G^{(N)}$ for (a) $N=3$; (b) $N=1$.
The fundamental domain $\Gamma_0$ is marked by a bold line.}
\end{figure}
%
%
%
%
%
%
%
%
%
%
For each integer $m\ge 1$ we define the finite group $\Z_m=\Z/(m\Z)$.
 The graph $\G^{(N)}$ is a union of edges $\G_\o$ and is given by
\[
\G^{(N)}=\cup_{\o\in \cZ} \G_\o,\ \ \o=(n,j,k)\in \cZ=\Z\ts \Z_3\ts \Z_N,\ \ 
\]
$$
\G_\o=\{r=r_\o^0+te_\o,\  t\in [0,1]\}, \qq e_\o=r_\o^1-r_\o^0,\qq r_\o^0,r_\o^1\in \R^3,\ \ \ |e_\o|=1,
$$
\[
\lb{ver1}
r_{\o_1}^0=r_{\o}^1=r_{\o_2}^0,\ \ where\ \ \o_1=(n+1,0,k),\ \o=(n,1,k),\ \o_2=(n,2,k),\  
\]
\[
\lb{ver2}
r_{\o_3}^1=r_{\o}^0=r_{\o_4}^1,\ \ where\ \ \o_3=(n,0,k),\ \ \o_4=(n,2,k-1),
\]
see Fig. \ref{fig1} and \ref{fig3}.
Each edge $\G_\o$ is a segment with length $|\G_\o|=1$.
We have the coordinate $r_\o=r_\o^0+te_\o$ and the local coordinate $t\in [0,1]$.   Introduce the vertex set of $\G^{(N)}$ by $V^{(N)}=\{v:v=r_\o^j, j=0,1, \o=(n,0,k)\in \cZ\}$.
For a function $f(x), x\in \G$ we define a function $f_\o=f|_{\G_\o}, \o\in \cZ$. We identify each function $f_\o$ on $\G_\o$ with a function on $[0,1]$ by using the local coordinate $t\in [0,1]$ in $r_\o=r_\o^0+t e_\o$. 
Our operator $\mH$ on the graph $\G^{(N)}$ acts
in the Hilbert space $L^2(\G^{(N)})=\sum_\o \os L^2(\G_\o)$. 
Then acting on the edge $\G_\o$, $\mH$ is the ordinary differential operator given by (see \cite{ARZ}, \cite{Ku})
\[
(\mH f)_\o=-\pa_{\o}^2f_\o(t)+q(t)f_\o(t),\qqq 
\pa_\o={d\/dt}-ia_\o,\qq a_\o(t)=(\mA(r_\o^0+te_\o),e_\o),
\]
where below in Sect. 2 we will show that 
\[
\label{aoI}
a_{n,0,k}=0,\qqq a=a_{n,1,k}=a_{n,2,k}={B\sqrt 3\/4}\cos {\pi\/2N},
\qqq all \qq (n,k)\in \Z\ts \Z_N.
\]
and $f_\o, f_\o'' \in L^2(\G_\o),\ \o\in \cZ$;   $q\in L^2(0,1)$ and $f\in \gD(\mH)$ satisfies the following

\begin{figure}
\centering
\psfrag{B}{$B$}
\includegraphics[width=.3\textwidth]{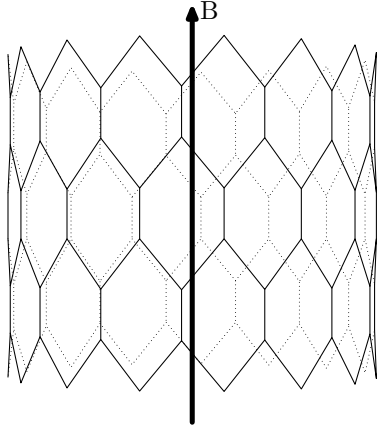}
\caption{The zigzag graph $\G$ subjected to the uniform magnetic field $B$.}
\label{fig3}
\end{figure}

\no {\bf The  Kirchhoff  Magnetic  Boundary Conditions:} {\it each $f\in \gD(\mH)$ is continuous on $\G^{(N)}$ and satisfies
\[
\lb{KirC}
-\pa_{\o_3}f_{\o_3}'(1)+\pa_{\o}f_{\o}'(0)-\pa_{\o_4}f_{\o_4}'(1)=0,\qqq
\pa_{\o_1}f_{\o_1}'(0)-\pa_{\o}f_{\o}'(1)+\pa_{\o_2}f_{\o_2}'(0)=0,
\]
$$
all \ 
\o_1=(n+1,0,k),\ \  \o=(n,1,k),\ \o_2=(n,2,k),\
\o_3=(n,0,k),\ \ \o_4=(n,2,k-1)\in \cZ,
$$
i.e., the sum of external magnetic derivatives of $f$ at each vertex  of $\G^{(N)}$ is equal to 0}.

In Theorem \ref{T1} we will show that the operator
$\mH$ is unitary equivalent to the operator $H=\sum_1^N\os H_k$.
The operator $H_k$ on the graph $\G^{(1)}$ acts
in the Hilbert space $L^2(\G^{(1)})$. 
In the case $N=1$ we will write  $\G_{n,j}=\G_{n,j,1}$ since  $k=1$
only. Thus $\G^{(1)}=\cup_{(n,j)\in \Z\ts \Z_3} \G_{n,j}$. Acting on the edge $\G_{n,j}$, $H_k$ is the ordinary differential operator given by 
\[
(H_k f)_{n,j}=-f_{n,j}''+q f_{n,j},\qqq f_{n,j}: [0,1]\to\C
\]
on the vector functions $f=(f_{n,j}), (n,j)\in\Z\ts \Z_3$, which
 satisfy the Kirchhoff conditions
\[
\lb{1K0}
f_{n,0}(1)=f_{n,1}(0)=e^{i a}s^k f_{n,2}(1),\quad
f_{n+1,0}(0)=e^{i a}f_{n,1}(1)=f_{n,2}(0),\qq s=e^{i{2\pi \/N}},
\]
\[
\lb{1K1}
-f'_{n,0}(1)+f'_{n,1}(0)-e^{i a}s^k f'_{n,2}(1)=0,\quad
f'_{n+1,0}(0)-e^{ia}f'_{n,1}(1)+f'_{n,2}(0)=0.
\]

Introduce the space $C(\G^{(N)})$  of complex continuous functions on $\G^{(N)}$ and the Sobolev space $W^2(\G^{(N)})=\rt\{f\in \C(\G^{(N)}): f\ \text{satisfies  conditions \er{1K0},\er{1K1}}, \ 
 f,f''\in L^2(\G^{(N)})\rt\}$.
If  $q=0$, then we denote our operator $H_{k}$ by $H_{0k}$. The operator $H_{0k}$ is self-adjoint on $\gD(H_{0k})=W^2(\G^{(1)})$[Ca1]. The operator $H_k$ is self-adjoint with $\gD(H_k)=\gD(H_{0k})$, see Sect.3.

In 1936 such operators were used to model aromatic molecules by Pauling \cite{Pa}. In 1953 ~Ruedenberg and Scherr \cite{RS}
described this model in details. We are going to consider in the framework of this model carbon nanotube, i.e. quasi-one-dimensional materials made of $sp^2$-hybridized carbon networks \cite{RS}.
Our graph $\G^{(N)}$ is a wrapped honeycomb lattice (see Fig \er{fig3}),
i.e., $\G^{(N)}$ is a single wall zigzag nanotube \cite{Ha}. 
For applications of our model spectral and analysis of these operators  see the liter. in \cite{ARZ}, \cite{Ku}.  

We reduce the spectral problem on the graph to
some matrix problem on $\R$. In order to describe this we 
define the fundamental subgraph $\G_0$ by
$
\G_{0}=\cup_{j=0}^2 \cup_{k=1}^N\G_{0,j,k},\ \ 
$
Thus $\cup_{k=1}^N\G_{0,0,k}$ contains only ''vertical`` edges $\G_{0,0,k}$;
$\cup_{k=1}^N\G_{0,1,k}$ and $\cup_{k=1}^N\G_{0,2,k}$ contains only edges $\G_{0,j,k}$with  positive and negative projections on the vector $(0,0,1)\in \R^3$, see Fig. \ref{fig1}. On $\G^{(N)}$, the group $\Z$ acts via
\[
\lb{grac}
p\circ  \G_{n,j,k}=\G_{n+p,j,k},\ \ n,p\in\Z,\ \ (n,j,k)\in \cZ. 
\]
Thus $\G_0$ is a fundamental domain associated with this group action.

There are two methods to study periodic differential operators.
The direct integral analysis usually used for partial
differential operators \cite{ReS} gives general information about the spectrum, but no detailed results.
The method of ordinary differential operators, based on the Floquet matrix analysis, gives detailed results, but even for 
the Schr\"odinger operator with periodic 2x2 matrix potentials
on the real line there are a lot of open problems \cite{BBK}.

We reduce the spectral problem on the graph to
some matrix problem on $\R$. For the operator $H_k$ we construct the fundamental solutions 
$
\vT_k(x,\l)
=(\vT_{k,\a}(x,\l))_{\a\in\cZ\ts \Z_3},$ and $ \F_k(x,\l)=(\F_{k,\a}(x,\l))_{\a\in\cZ\ts \Z_3}$, $(x,\l)\in \R\ts \C
$ 
which  satisfy
\[
\lb{eqf}
-f_{n,j}''+qf_{n,j}=\l f_{n,j},\ \ 
\ \ \text{the Kirchhoff Boundary Conditions \er{1K0},\er{1K1}}, 
\]
\[
\lb{eqf0}
\vT_{k,\a}(0,\l)=\F_{k,\a}'(0,\l)=1,\qq 
\vT_{k,\a}'(0,\l)=\F_{k,\a}(0,\l)=0,\qqq \a=(0,0).
\]
We introduce the monodromy matrix
\[
\cM_k(\l)=\ma \vT_{k,\o}(0,\l) & \F_{k,\o}(0,\l)\\
\vT_{k,\o}'(0,\l) & \F_{k,\o}'(0,\l)\am,\qqq \o=(1,0).
\]

We introduce the monodromy matrix, similar to the case
of Schr\"odinger operators with periodic matrix potentials
on the real line, see \cite{YS}.
After this, roughly speaking, using the approach from \cite{BK},\cite{BBK},\cite{CK} we introduce the Laypunov functions and study   
the properties of these functions, similar to the 
case of  the Schr\"odinger operator with a periodic matrix
potential on the real line. 
This is a crucial point of our analysis.
 Here we essentially use the results and techniques
from the papers \cite{BK}, \cite{BBK},\cite{CK}.
The recent papers \cite{BBK},\cite{CK} are devoted to the 
the Schr\"odinger operator with a periodic matrix
potential (a standard case) on the real line.
Remark that Carlson \cite{Ca} studied the monodromy 
operator to analyze the Schr\"odinger operator 
on  a product of graphs.

Recall the needed properties of 
the equation $-f''+q(x)f=\l f$ on the real line
with a periodic potential $q(x+1)=q(x),x\in \R$. In this case we 
introduce the fundamental solutions $\vt(x,\l)$ 
and $\vp(x,\l),x\in \R$ satisfying $\vt(0,\l)=\vp'(0,\l)=1, \vt'(0,\l)=\vp(0,\l)=0$.
The corresponding monodromy matrix $\cM(\l)$
and the Lyapunov function $\D$  are given by
\[
\cM(\l)=\ma\vt(1,\l) & \vp(1,\l) \\ 
\vt'(1,\l) & \vp'(1,\l)\am,\ \ \ \ \ 
\D(\l)={1\/2}(\vp'(1,\l)+\vt(1,\l)),\ \ \l\in\C.
\]
Let $\m_n, n\ge 1,$ be the Dirichlet spectrum of 
the problem $-y''+qy=\l y, y(0)=y(1)=0$ on the unit interval
$[0,1]$ and define the set $\s_D=\{\m_n, n\ge 1\}$. Recall that $\s_D=\{\l\in\C: \vp(1,\l)=0\}$. Let $\n_n, n\geq 0,$ be the Neumann
 spectrum of equation $-y''+qy$  for the boundary
 condition $y'(0)=y'(1)=0$. Define 
\[
F=2\D^2+{\vp(1,\cdot)\vt'(1,\cdot)\/4}-1,\ \qqq
c_k=\cos(a+{\pi k\/N}),\qqq s_k=\sin(a+{\pi k\/N}),\qq
\]
$k=0,1,..,N$. We formulate our first result about the fundamental
solutions $\vT_k, \F_k$.

\begin{theorem}
\label{T1}
i) The operator $\mH$ is unitary equivalent to
the operator $H=\sum_1^N\os H_k$.

\no ii) Let $c_k=\cos(a+{\pi k\/N})\ne 0$ for some $(k,a)\in\Z_N\ts\R$.
Then for any $\l\in\C\sm\s_D$  there exist unique fundamental solutions $\vT_k, \F_k$ of the system \er{eqf} with conditions \er{eqf0} and  each function $\vT_k(x,\l),\F_k(x,\l),x\in \G^{(1)}$ is meromorphic in $\l\in\C\sm\s_D$.  Moreover, each matrix $\cM_k(\l)$ satisfies
\[
\lb{T1-1}
\cM_k=\cR^{-1}\cT_k \cR \cM ,\quad
\cT_k={s^{-{k\/2}}\/2c_k}\ma 2\D & 1 \\ 4\D^2-4c_k^2 & 2\D \am,\qq 
\cR=\ma 1 & 0\\0 & \vp(1,\cdot)\am,
\]
\[
\label{T1-2}
\det \cM_k={s}^{-k},\qqq
\Tr\cM_k={2s^{-{k\/2}}(F+s_k^2)\/c_k}
\]
and the function $D_k(\t,\l)=\det(\cM_k(\l)-\t I_{2})$ is entire
with respect to $\l, \t\in \C$.

\end{theorem}

Remark that in contrast to the Schr\"odinger operator with periodic matrix potentials on the real line (see \cite{YS} or \cite{CK}),
the monodromy matrix $\cM_k$ has poles at the points $\l\in \s_D$,
which are eigenvalues of $H_k$, see Theorem \ref{T2}.
However, $\cM_k$ is similar to the entire matrix-valued function.
Define the subspace $\cH_k(\l)=\{\p\in \gD(H_k): H_k\p=\l \p\}$
for $\l\in \s(H_k), k\in \Z_N$. If some $\l_0\in \R$ is 
an eigenvalue of $H_k$
with infinite multiplicity, then we say that $\{\l_0\}$ is a flat band. In Theorem \ref{T2} and \ref{T3} we describe all flat bands.

\begin{theorem}\lb{T2} 
Let $(\l,a,k)\in\s_D\ts\R\ts \Z_N$. Then

\no i) Every eigenfunction from $\cH_k(\l)$ vanishes at all vertexes of $\G^{(1)}$.

\no ii) Let the function $\p^{(0)}$ be given by :

if $\e=1-e^{2ia}s^k\vp'(1,\l)^2\ne 0$, then
\begin{multline}
\lb{T2-1}
\p^{(0)}_{n,j}=0, \text{for all}\  n\ne 0,-1, \ j\in \Z_3,\qqq
and \qq
\p^{(0)}_{0,0}=\e\vp_t, \ \p^{(0)}_{0,1}=\vp_1'\vp_t, \ 
\p^{(0)}_{0,2}=e^{ia}{\vp_1'}^2\vp_t, \\
\p^{(0)}_{-1,0}=0,\ 
\p^{(0)}_{-1,1}=-e^{ia}s^k{\vp_1'}\vp_t, \ \p^{(0)}_{-1,2}=-\vp_t, \ t\in [0,1],
\end{multline}

if $\e=0$, then
\[
\lb{T2-2}
\p^{(0)}_{0,0}=0,\qq
\p^{(0)}_{0,1}(t)=\vp_t,\qq
\p^{(0)}_{0,2}(t)=e^{ia}\vp_1'\vp_t,\qq t\in[0,1], 
\p^{(0)}_{n,j}=0,\ all \ n\ne 0, j\in \Z_3.
\]
Then  each $\p^{(n)}=(\p^{(0)}_{n-m,j})_{(m,j)\in \Z\ts \Z_3}\in \cH_k(\l), n\in \Z$  and each $f\in \cH_k(\l)$ has the form
\[
\lb{T2-3}
f=\sum_{n\in\Z}\wh f_n\p^{(n)},\qqq 
\wh f_{n}=\ca\e^{-1} f_{n,0}'(0) & if \ \e\ne 0\\
f_{n,1}'(0) & if \ \e= 0 \ac , 
\qqq (\wh f_n)_{-\iy}^\iy\in \ell^2
\]
and the mapping $f\to \{\wh f_{n}\}$ is a linear isomorphism between $\cH_k(\l)$ and $\ell^2$.

\end{theorem}

\begin{figure}
\tiny
\centering
\noindent
(a)\quad\quad\includegraphics[angle=90,width=.4\textwidth]{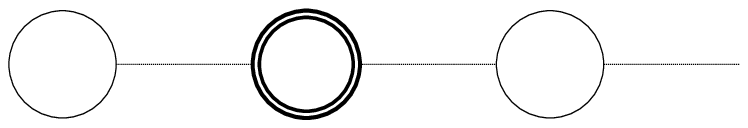}
(b)\quad\quad\includegraphics[angle=90,width=.4\textwidth]{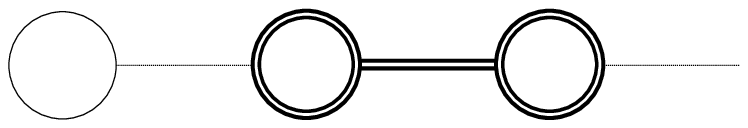}
\caption{The supports of the eigenfunction $\p^{(0)}$ corresponding
to an eigenvalue contained in $\s_D$: (a) $\e=1-e^{2ia}s^k\vp'(1,\l)^2=0$; (b) $\e\neq 0$.
The support of the eigenfunction $\p^{(0)}$ corresponding to
an eigenvalue belonging to $\s_{AP}$ coincides with
one shown on the figure (b).
}
\lb{fig2}
\end{figure}

For a self-adjoint operator $H$ we define the set
$$
\s_{\iy}(H)=\{\l: \ \l\in\s_{pp}(H)\ \text{is of 
infinite multiplicity} \}.
$$
We define the {\bf singular magnetic field} by $\wt a_{k,m}={\pi\/2}-{\pi k\/N} +\pi m, (k,m)\in Z_N\ts \Z$ and  {\bf singular energy} $\l\in\s_{AP}
=\{\l\in \R: F(\l)=-1\} $. Note that $\s_D\cap \s_{AP}=\es$, see Sect.4 and Theorem \ref{T04}. Let $\wt\l_n, n\ge 1$ be the zeros of $F(\l)=-1$.
Due to Theorem \ref{T04}, all $\wt\l_n$ are real and have asymptotics $\wt\l_n=\pi^2 n^2+O(n)$ as $n\to \iy$. We describe the operator $H_k(a)$ for the case of {\bf singular magnetic field}.

\begin{theorem}\lb{T3}
Let $c_k=\cos(a+{\pi k\/N})=0$ for some $(k,a)\in\Z_N\ts\R$. Then 
\[
\lb{T3-1}
\s(H_k(a))=\s_{\iy}(H_k(a))=\s_D\cup \s_{AP},\qq where \qq \s_{AP}
=\{\l\in \R: F(\l)=-1\} .
\]
If in addition for some $\l\in \s_{AP}$ a function $\p^{(0)}$ is given by
\begin{multline}
\lb{T3-2}
\p^{(0)}_{0,0}(t)=\vp_t C_1+\vt_t C_2,\qqq 
\p^{(0)}_{0,1}(t)=\vt_t-\vp_t{\vt_1\/\vp_1},\qqq 
\p^{(0)}_{0,2}(t)=-{\vp_t\/\vp_1}e^{ia},\\
\p^{(0)}_{-1,0}=0,\qq \f_{-1,1}(t)={\vp_t\/\vp_1}e^{-ia}C_1,\qq 
\p^{(0)}_{-1,2}(t)=(\vt_t-\vp_t{\vt_1\/\vp_1})C_1,
\qq t\in[0,1],
\end{multline}
\[
\lb{T3-3}
\p^{(0)}_{n,j}=0, any \ n\ne -1, 0, j=1,2,\qq
C_1=\vp_1'+2\D,\qq  C_2=-\vt_1'-{2\vt_1\D\/\vp_1},
\]
then each $\p^{(n)}=(\p^{(0)}_{n-m,j})_{(m,j)\in \Z\ts \Z_3}\in \cH_k(\l), n\in \Z$ and each $f\in \cH_k(\l)$ has the form
\[
\lb{T3-4}
f=\sum_{n\in\Z}\wh f_n\p^{(n)},\qqq 
\wh f_{n}={f_{n,1}(0)\/C_2}, \qqq (\wh f_n)_{-\iy}^\iy\in \ell^2
\]
and the mapping $f\to \{\wh f_{n}\}$ 
is a linear isomorphism between $\cH_k(\l)$ and $\ell^2$.
\end{theorem}

Below we will sometimes write $\cM_k(\l,a,q), F_k(\l,a,q),\dots$, instead of $\cM_k(\l), F_k(\l),\dots$, when several magnetic fields and potentials are being dealt with. Let $\t_{k,\pm}$ be the eigenvalues of $\cM_k, k\in \Z_N$. Using \er{T1-2} we deduce that
\[
\lb{malk}
\t_{k,-}\t_{k,+}=s^{-k}, \qq \t_{k,-}+\t_{k,+}=\Tr 
\cM_{k}={2s^{-{k\/2}}\/c_k}(F+s_k^2).
\]
If $k=0$, then we can introduce the standard entire Lyapunov function 
$F_0$ by
\[
\lb{DeL0}
F_0={\Tr \cM_0\/2}={1\/2}(\t_{0,+}+{1\/\t_{0,+}})={F+s_0^2\/c_0}.
\]
The asymptotics of the fundamental solutions $\vt, \vp,$ imply
(see [KL])
$$ 
\lb{asD0}
F(\l,q)=\D_0^0(\l)+{O(e^{|\Im\sqrt{\l}|})\/\sqrt{\l}},\qq
\D_0^0(\l)={9\cos 2\sqrt{\l}-1\/8}\qq as \qq |\l|\to \iy,  
$$
 where $\D_0^0$ is $F$ at $q=0$.
In particular, if $q=0$, then we get $F_0(\l,a,0)={\D_0^0+s_0^2\/c_0}$.

If $k\neq 0$, then we define the Lyapunov functions 
$F_{k,\pm}={1\/2}(\t_{k,\pm}+{1\/\t_{k,\pm}})$, (see \cite{BBK},\cite{CK}). 
But in this case the Lyapunov functions are not entire, in general (see \cite{BBK},\cite{CK}). Below we prove the following identities
\[
\lb{DeLk}
F_{k,\pm}={1\/2}(\t_{k,\pm}+{1\/\t_{k,\pm}})=T_k\pm\sqrt{R_k},\qq T_k={c_{0k}\/c_{k}}(F+s_k^2),\qqq
R_k={s_{0k}^2\/c_{k}^2}\rt(c_k^2-(F+s_k^2)^2\rt),
\]
where $s_{0,k}=\sin{\pi k\/N},c_{0,k}=\cos{\pi k\/N}$.
Introduce the two sheeted Riemann surface $\gR_k$
(of infinite genus) defined by $\sqrt {R_k}$.
The functions $F_{k,\pm},k\in \ol {N-1}$ 
are the branches of $F_k=T_k+\sqrt{R_k}$ on the Riemann surface $\gR_k$, where the set $\ol N=\{1,2,..,N\}$.

\begin{theorem}
\lb{T4} 

Let $c_k=\cos(a+{\pi k\/N})\ne 0$ for some $(k,a)\in\Z_N\ts\R$. Then

\no i) The Lyapunov functions $F_{k,\pm}, k\ne 0$ satisfy \er{DeLk}. 
 
\no ii) The  following identities hold:
\[
\lb{T4-1}
\s(H_k)=\s_{\iy}(H_k)\cup\s_{ac}(H_k),\qq 
\s_{\iy}(H_k)=\s_D,\qq \s_{ac}(H_k)=\{\l\in\R: F_k(\l)\in [-1,1]\}.
\] 
iii) Let some $\l\in \R$ be not a branch point of $F_k$ and let
$F_k(\l)\in (-1,1)$. Then $F_k'(\l)\neq 0$. 

\end{theorem}

{\bf Remark.} 1) If we know $F_0$, then we determine all
$F_k, \r_k,k=\ol {N-1}$ by \er{DeLk}. 2) If we know $\r_k$ for some $k\in\ol {N-1}$, then we determine all
$F_k, \r_k,k\in\ol {N}$ by \er{DeLk}.

Consider $H_0$. Using the definition \er{DeL0} we obtain  $F_0(\cdot,-a)=F_0(\cdot,a)=-F_0(\cdot,a+\pi)$
for $\cos a\ne 0$ and
$
F_0(\cdot,a+{\pi\/2})=-F_0(\cdot,{\pi\/2}-a)$ for 
$a\in \C\sm \pi \Z$. Due to this symmetry it is sufficient to
consider the case $a\in [0,{\pi\/2})$. The value
$\wt a_{0,m}={\pi\/2}+\pi m$ is a singular magnetic field for $H_0$.
The important case $a=0$
for odd $N$ was studied in \cite{KL}, the case of $N$ even is studied in Sect. 4. The results for both these cases are formulated in  Sect. 4.

\begin{figure}
\label{fig4}
\centering
\tiny
\psfrag{54}[c][c]{$-\frac{5}{4}$}
\psfrag{-1}[c][c]{$-1$}
\psfrag{1}[c][c]{$1$}
\psfrag{0}[c][c]{$0$}
\psfrag{Pi/2}[t][b]{$\l_{0,1}$}
\psfrag{Pi}[t][b]{$\l_{0,2}$}
\psfrag{(3*Pi)/2}[t][b]{$\l_{0,3}$}
\psfrag{2*Pi}[t][b]{$\l_{0,4}$}
\psfrag{l}{$\l$}
\noindent
\begin{tabular}{cc}
(a)\quad\quad\includegraphics[width=.4\textwidth]{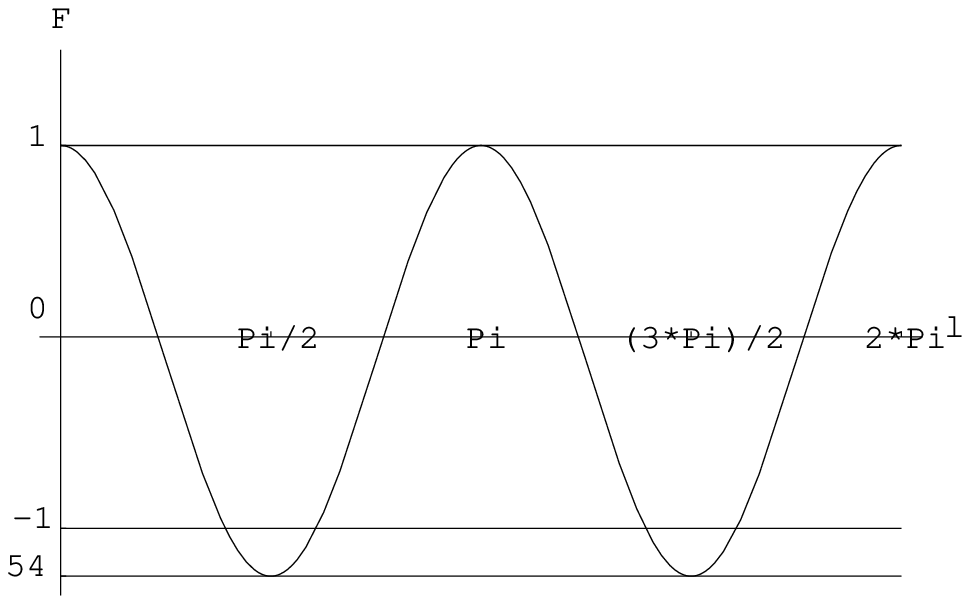}&
(b)\quad\quad\includegraphics[width=.4\textwidth]{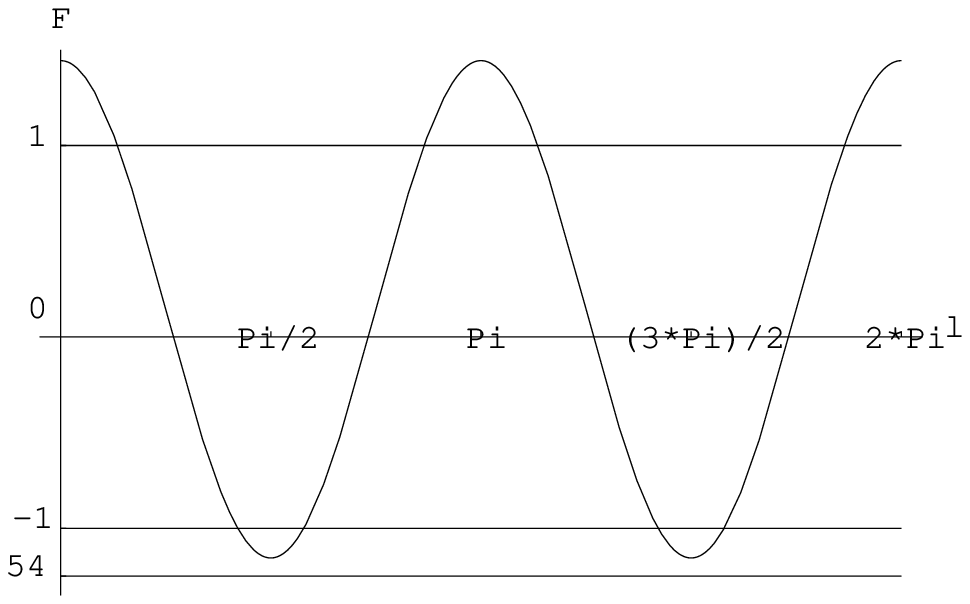}\\
(c)\quad\quad\includegraphics[width=.4\textwidth]{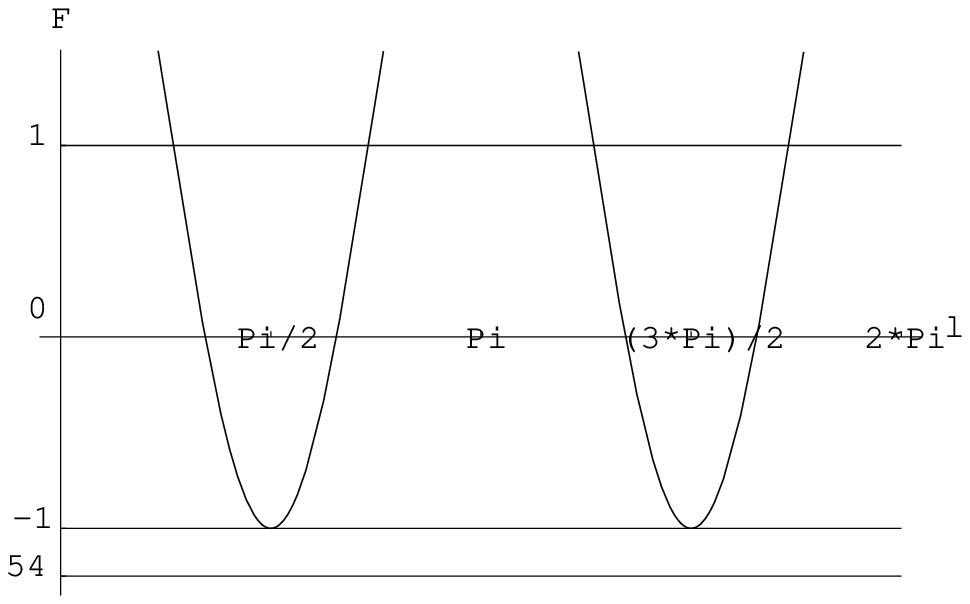}&
(d)\quad\quad\includegraphics[width=.4\textwidth]{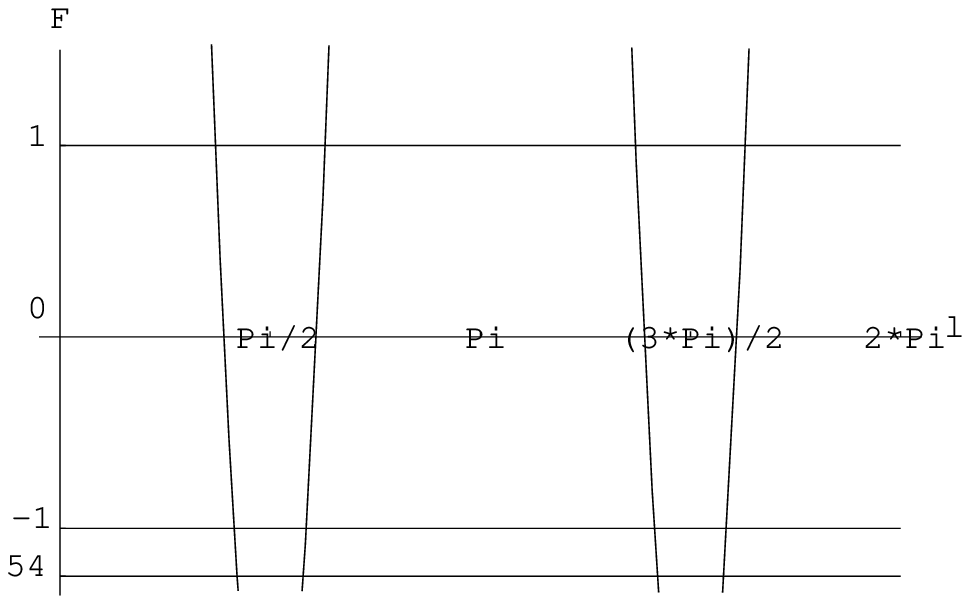}
\end{tabular}
\caption{The function $F_0(\l,a)$ for (a) $a=0$; (b) $a=\frac{\pi}{6}$;
(c) $a=\frac{\pi}{3}$; (d) $a\uparrow \frac{\pi}{2}$.}
\end{figure}


Let $a\in (0,{\pi\/2})$.
Let $D_0(\cdot,\t)=\det (\cM_0-\t I_2),\t\in \C$ and let $\l_{0,n}^\pm=\l_{0,n}^\pm(a)$ be the zeros of $D_0(\l,1)D_0(\l,-1)$, where $I_N,N\ge 1$ is the identity $N\ts
N$ matrix. The zeros of $D_0(\l,1)$ ( and $D_0(\l,-1)$) counted
with multiplicity are the periodic (anti-periodic) eigenvalues
for the equation $-y''+qy=\l y$ on $\G^{(1)}$  with periodic (anti-periodic) boundary conditions. Below we will show that $\l_{0,n}^\pm=\l_{0,n}^\pm(a)$ satisfy for $a\in (0,{\pi\/2})$:
\[
\lb{1epa0}
F_0(\l_{0,n}^\pm)=(-1)^n,\ \ \ 
\l_{0,0}^+<\l_{0,1}^-< \l_{0,1}^+<\l_{0,2}^-< \l_{0,2}^+<\l_{0,3}^-< \l_{0,3}^+<\l_{0,4}^-< \l_{0,4}^+<....
\]
\[
\lb{Tas-4}
\l_{0,n}^\pm=({\pi n\/2}\pm \f_{n})^2+q_0+{o(1)\/n},\ \ 
\qq \f_{2n}=f_0,\ \f_{2n-1}={\pi \/2}-f_1\qq  as \qq n\to \iy,
\]
where $f_{s}={1\/2}\arccos {1+(-1)^s8c_{0}-8s_0^2\/9}\in [0,{\pi\/2}],
s=0,1$. If $a=0$, then $\l_{0,2n}^-\le \l_{0,2n}^+,n\ge 1$, see Sect.4. Introduce the gaps $\g_{0,n}(a)$ and the spectral bands $\s_{0,n}(a)$ of $H_0(a)$ by
$$
\g_{0,n}(a)=(\l_{0,n}^-(a),\l_{0,n}^+(a)),\qqq \s_{0,n}(a)=[\l_{0,n-1}^+(a),\l_{0,n}^-(a)], \qq n\ge 1.
$$
The following theorem describes the basic properties of $H_0$.

\begin{theorem}\label{T5}
Let $a\in (0,{\pi\/2})$. Then the function $F_0(\l,a)=(F(\l)+s_0^2)/c_0$ satisfies

\no i) The function ${\pa\/\pa \l} F_0(\l,a)$ has only real simple zeros
$\l_{0,n}=\l_{0,n}(a),n\ge 1$,
which are separated by the simple zeros $\e_{0,n}=\e_{0,n}(a)$ of $F_0(\l,a)$:
$\e_{0,1}<\l_{0,1}<\e_{0,2}<\l_{0,2}<\e_{0,3}<...$ and satisfy
for all $n\ge 1$:
\[
\lb{T5-1}
\n_n,\m_n\in \g_{0,2n},\qq 
F_0(\l_{0,2n},a)\ge {1+s_0^2\/c_0}>1,\qqq   F_0(\l_{0,2n-1},a)\le -{1+4c_0^2\/4c_0}\le -1.
\]
Moreover, if $c_0\ne {1\/2}$, then  $F_0(\l_{0,2n-1},a)\le -{1+4c_0^2\/4c_0}< -1$. The periodic and anti-periodic eigenvalues
satisfy \er{1epa0}.

\no ii) Each even gap $\g_{0,2n}(a)=(\l_{0,2n}^-,\l_{0,2n}^+), n\ge 1$ is open and satisfies
\[
\lb{T5-2}
 \g_{0,2n}(a)\ss \g_{0,2n}(a_1), \qqq 0<a<a_1<{\pi\/2}.
\]
Odd gaps $\g_{0,m}(a)=(\l_{0,m}^-(a),\l_{0,m}^+(a)),m=2n-1, n\ge 1$ satisfy
\[
\lb{T5-3}
 \g_{0,m}(a)\supset \g_{0,m}(a_1), \qq 0\le a< a_1\le{\pi\/3}\ \ and \ \
 \g_{0,m}(a)\ss \g_{0,m}(a_1), \qq {\pi\/3}\le a<a_1<{\pi\/2}.
\]
Moreover, if $a={\pi\/3}$ and $q\in L_{even}^2(0,1)$, then all $\g_{0,2n-1}(a)=\es, n\ge 1$.

\no iii) The asymptotics \er{Tas-4} hold true.

\no iv) Let $a\uparrow {\pi\/2}$. Then  each band       $\s_{0,n}(a)\to\wt\s_n=\{\wt\l_n\}$ (recall $F(\wt\l_n)=-1$) and $\wt\s_n$ is a flat
 band for $H_0({\pi\/2})$ and 
 \[
\lb{T5-4}
\l_{0,n}^-(a), \ \l_{0,n-1}^+(a)=\wt\l_n+{c_0 +O(c_0 ^2)\/ F'(\wt\l_n)}\qqq
 {\rm as} \qq a\uparrow{\pi\/2}.
\]
\end{theorem}

{\bf Remark.} 1) If $a={\pi\/2}$, then by Theorem \ref{T3},
$\s(H_0({\pi\/2}))=\s_\iy(H_0({\pi\/2}))=\s_{AP}\cup \s_D$.
By Theorem \ref{T4}, \ref{T5}, $\s(H_0(a))=\s_{ac}(H_0(a))\cup \s_\iy(H_0(a)), \s_\iy(H_0(a))=\s_D$ for all $a\in (0,{\pi\/2})$.
The asymptotics \er{T5-2} give the asymptotics of shrinking bands,
as $a\uparrow{\pi\/2}$. Moreover, if ${\pi\/2}-a>0$ is increasing
then the spectral bands are increasing too.  

2) If $a\in (0,{\pi\/2})$,
then all even gaps $\g_{0,2n}$ are open. If $a\uparrow {\pi\/2}$,
then each band shrinks to the flat band.

3) If firstly $a={\pi\/2}$, then by Theorem \ref{T3}, all bands
are flat.  If secondly, $a<{\pi\/2}$, then by iv), each flat band $\wt\s_n$ becomes ordinary band from $\s_{ac}(H(a))$.

Consider the operator $H_k,k\ne 0$ for the case $c_k\ne 0$. Let $D_k(\cdot,\t)=\det (\cM_k-\t I_2),\t\in \C, k\in \ol{N-1}$.
Let $\l_{k,n}^\pm=\l_{k,n}^\pm(a)$ and $\m_{k,n}^\pm=\m_{k,n}^\pm(a),n\ge 0$ be the zeros of $D_k(\l,1)$ and $D_k(\l,-1)$. 
In Theorem \ref{T6}  we will show that the periodic  eigenvalues $\l_{k,n}^\pm$ and the anti-periodic eigenvalues $\m_{k,n}^\pm$ satisfy  the equations $F_k(\l_{k,n}^\pm)=1,\ F_k(\m_{k,n}^\pm)=-1$ and
\[
\lb{epak1}
F(\l_{k,n}^\pm)=c_{0k}c_k-s_k^2\in (-{5\/4},1), \qq  
F(\m_{k,n}^\pm)=-c_{0k}c_k-s_k^2\in (-{5\/4},1),
\]
and a labeling is given by: each $\l_{k,n}^\pm, \m_{k,n}^\pm$ is simple and
\[
\lb{epak2}
\l_{k,0}^+<\l_{k,1}^-<\l_{k,1}^+
<\l_{k,2}^-<\l_{k,2}^+...., \qqq
\m_{k,0}^+<\m_{k,1}^-<\m_{k,1}^+
<\m_{k,2}^-<\m_{k,2}^+...., 
\]
\[
\lb{epak3} 
\l_{k,n}^{\pm}=(\pi n\pm \f_{k,0})^2+q_0+{o(1)\/n},\ \ \qq 
\m_{k,n}^{\pm}=(\pi n\pm \f_{k,1})^2+q_0+{o(1)\/n},\qq as\qq n\to \iy,
\] 
$$
\f_{k,s}={1\/2}\arccos {1+8((-1)^sc_{0k}c_k-s_k^2)\/9}\in [0,{\pi\/2}] \qq
\qq s=0,1.
$$
Let $r_{k,n}^\pm, k\in\ol{N-1}, n\ge 0$ be the zeros of $R_k$, which defined by \er{DeLk}.
A zero of $R_k, k\in \ol{N-1}$ is called a {\bf resonance} of $H_k$. 
Below we will show that these resonances satisfy the equations
\[
\lb{eqr}
F(r_{k,2n}^\pm )=|c_k|-s_k^2, \  \qq 
F(r_{k,2n+1}^\pm )=-|c_k|-s_k^2, \ \qq
\ \ \ k\in \ol{N-1},\ n\ge 0,
\]
they are real and labeling is given by
\[
\lb{esr}
r_{k,0}^+<
r_{k,1}^-<r_{k,1}^+< r_{k,2}^-<r_{k,2}^+ <..,\ \ \ \ |c_k|\neq {1\/2},
|c_k|\neq 1,
\]
\[
\lb{esr0}
r_{k,0}^+<
r_{k,1}^-<r_{k,1}^+< r_{k,2}^-\le r_{k,2}^+ <
r_{k,3}^-<r_{k,3}^+< r_{k,4}^-\le r_{k,4}^+ <..,\ \ \ \ |c_k|=1,
\]
\[
\lb{esrI}
r_{k,0}^+<
r_{k,1}^-\le r_{k,1}^+< r_{k,2}^-<r_{k,2}^+ <r_{k,3}^-\le r_{k,3}^+<..,
\qqq |c_k|={1\/2}.
\]
The resonances $r_{k,n}^{\pm},$ satisfy 
\[
\lb{Tas-r}
r_{k,n}^\pm=({\pi n\/2}\pm \f_{n})^2+q_0+o(n^{-1}),\ \ 
\f_{2n}=f_{0}, \f_{2n-1}={\pi\/2}-f_{1},\ 
\qq f_{s}={\arccos Y_s\/2} \in [0,{\pi\/2}], 
\]
$Y_s={1+(-1)^s8|c_{k}|-8s_k^2\/9}, s=0,1$ as $n\to \iy$.
We describe the spectral properties of $H_k$ in terms of 
the Lyapunov functions.

\begin{theorem}
\lb{T6} 
i) Let $\gS_k(a)=\{\l\in \R: (F(\l)+s_k^2)^2\le c_k^2\}$ 
for some $(k,a)\in \Z_N\ts \R$.
 Then the following identities hold
\[
\lb{T6-1}
\s(H_k(a))=\gS_k(a)\cup \s_D
=\s(H_{k-1}(a+{\pi\/N})),\qq
\gS_k(a)=\ca\s_{ac}(H_k(a)) &\ if\qq c_k\ne 0\\
\s_{AP} & if \qq c_k=0\ac.
\]
\no ii) Let $c_k=\cos(a+{\pi k\/N})\ne 0$  for some $(k,a)\in\ol{N-1}\ts\R$. Then the resonances  $r_{k,n}^\pm, n\ge 0$ satisfy Eq. \er{eqr} and the estimates \er{esr}-\er{Tas-r}.
Moreover, the spectral bands $\s_{k,n}$ and the gaps $\g_{k,n},n\ge 1$
for the operator $H_k$ are given by
\[
\lb{T6-2}
\s_{ac}(H_k(a))=\cup_{n=1}^\iy \s_{k,n}(a),\qq
\s_{k,n}(a)=[r_{k,n-1}^+(a),r_{k,n}^-(a)],\qq 
\g_{k,n}(a)=(r_{k,n}^-(a),r_{k,n}^+(a)).
\]
\no  iii)  Let  $c_k\ne 0$  for some $(k,a)\in \ol{N-1}\ts\R$.
The periodic $\l_{k,n}^\pm$ and anti-periodic eigenvalues $\m_{k,n}^\pm,n\ge 0$ satisfy  \er{epak1}- \er{epak3} and $\l_{k,n}^\pm, \m_{k,n}^\pm\in \s_{k,n}(a)$.

\no iii) Let $a\to \wt a_{k,m}={\pi\/2}-{\pi k\/N}+\pi m$ for some $(k,m)\in \ol{N-1}\ts \Z$. Then  each band $\s_{k,n}(a)\to\wt\s_n=\{\wt\l_n\}$ and $\wt\s_n$ is a flat band for $H_k(\wt a_{k,m})$ and
 \[
\lb{T6-3}
r_{k,n-1}^+(a), r_{k,n}^-(a)=\wt\l_n+{c_0 +O(c_0 ^2)\/ F'(\wt\l_n)}\qqq
 {\rm as} \qq a\to \wt a_{k,m}.
\]
\end{theorem}

\begin{figure}
\centering
\tiny
\psfrag{54}[c][c]{$-\frac{5}{4}$}
\psfrag{-1}[c][c]{$-1$}
\psfrag{1}[c][c]{$1$}
\psfrag{0}[c][c]{$0$}
\psfrag{Pi/2}[t][b]{$\frac{\pi^2}{4}$}
\psfrag{Pi}[t][b]{$\pi^2$}
\psfrag{(3*Pi)/2}[t][b]{$\frac{9\pi^2}{4}$}
\psfrag{2*Pi}[t][b]{$4\pi^2$}
\psfrag{l}[c][c]{$\l$}
\psfrag{F}[c][c]{}
\noindent
\begin{tabular}{cc}
(a)\quad\quad\includegraphics[width=.4\textwidth]{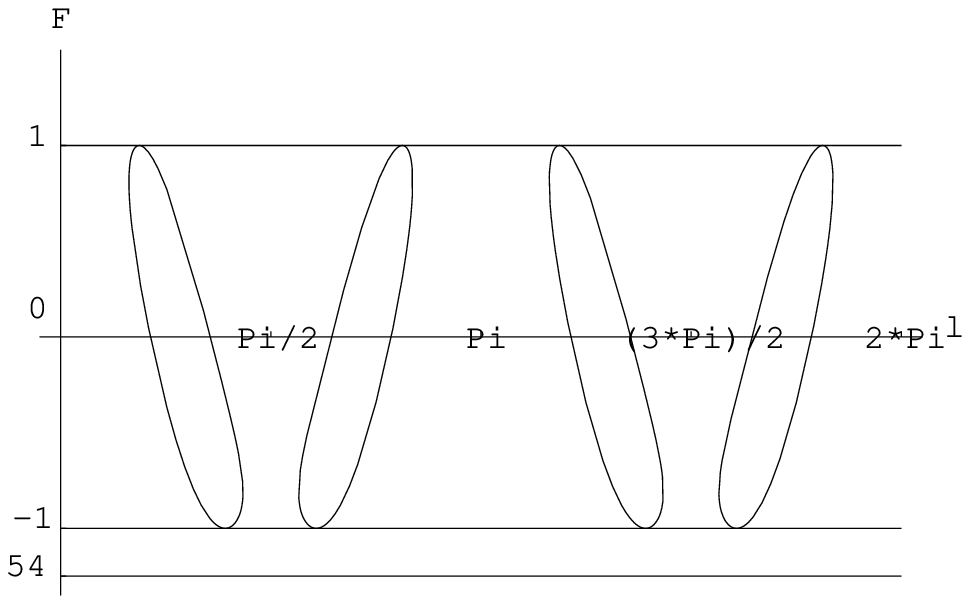}&
(b)\quad\quad\includegraphics[width=.4\textwidth]{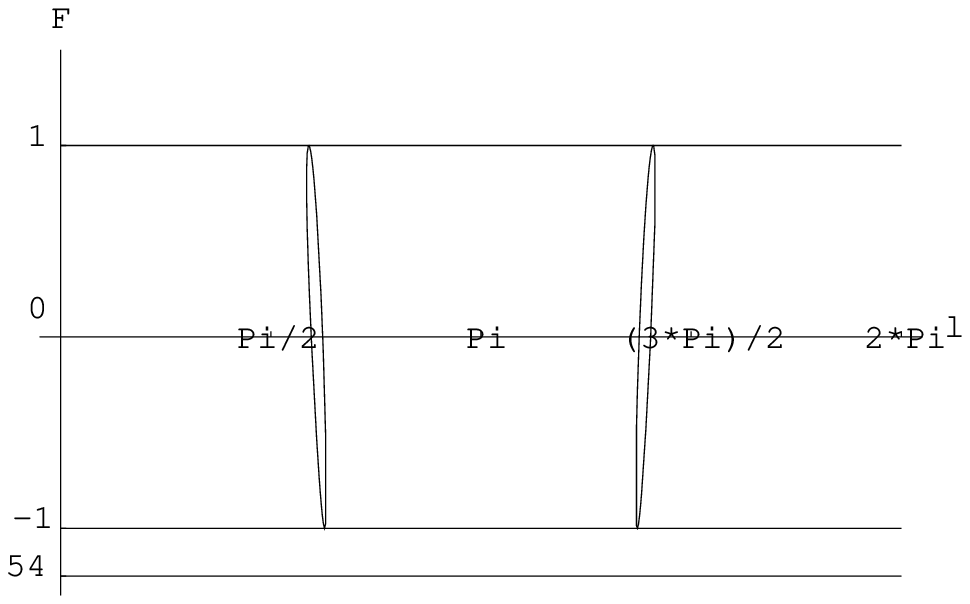}\\
(c)\quad\quad\includegraphics[width=.4\textwidth]{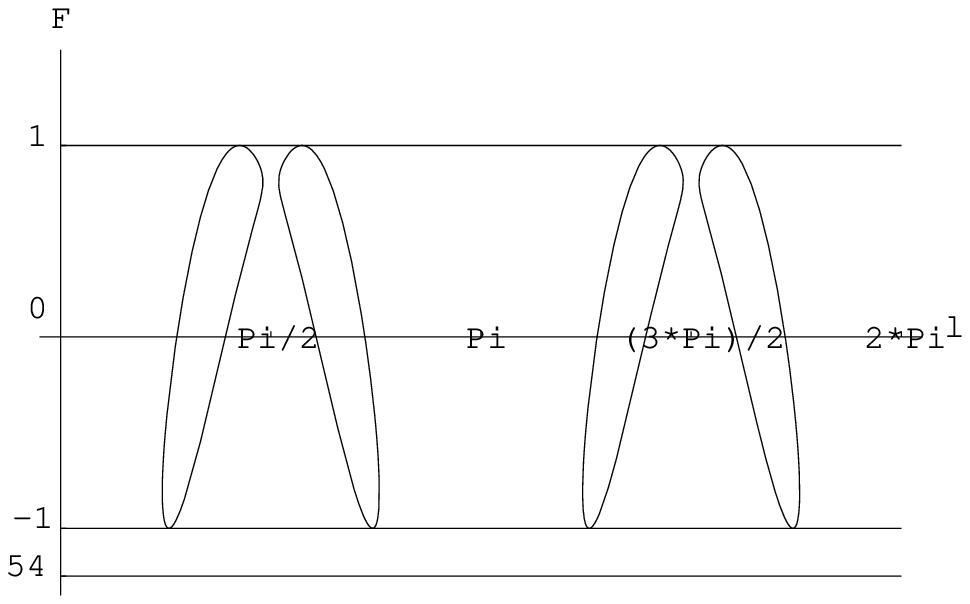}
\end{tabular}
\caption{The function $F_{1,\pm}(\l,a)$ for (a) $a=0$; (b) $a=\frac{\pi}{3}$;
(c) $a=\frac{\pi}{2}$.}
\label{fig5}
\end{figure}

Note that in the matrix case (see \cite{BBK}) 
the resonances in general have non-zero imaginary part.
The fact that all our resonances are  real is a peculiarity of the high
symmetry of a zigzag periodic graph.

The spectrum of the operator $H=H(a),a\in \R$ consists of an absolutely continuous part (spectral band $S_n(a),n\ge 1$ separated by gaps $G_n(a),n\ge 1$) plus an infinite number of eigenvalues  with infinite multiplicity, and we will show the following identities
\[
\lb{decsp1}
\s(H(a))=\s_{ac}(H(a))\cup \s_{\iy}(H(a)),
\]
\[
\lb{decsp2}
\s_{ac}(H(a))=\cup_{n\ge 1} S_n(a),\qq
S_n(a)=[E_{n-1}^+(a),E_n^-(a)],\qq 
G_n(a)=(E_{n}^-(a),E_n^+(a)),
\]
where $S_n(a),G_n(a)$ are given by \er{T7-5}, \er{T7-6}.

\begin{theorem}
\lb{T7} 
i) The following identities hold
\[
\lb{T7-1}
\s(H(a+{\pi\/N}))=\s(H(a)),\qqq a\in\R.
\]
\no ii)  Let $a\in (0,{\pi\/N})$. Then all gaps $G_n(a)=(E_{n}^-(a),E_{n}^+(a))$ for even $n\ge 2$ satisfy
\[
\lb{T7-2}
G_n(0)\ss G_n(a)\ss G_n(a'),\qqq if \qq 0<a<a'<{\pi\/N},
\]
\[
\lb{T7-3}
F(E_{n}^\pm(a))=c_+-1+c_+^2,\qq
c_+=\cos a_+,\qqq a_+=\max \{a,{\pi\/N}-a\}.
\]
In particular, if $a_+=a$, then $E_{n}^\pm(a)$ are periodic eigenvalues.
If $a_+={\pi\/N}-a$, then $E_{n}^\pm(a)$ are resonances. 
Moreover,
\[
\lb{T7-4}
E_{n}^\pm(a)=({\pi n\/2}\pm \f_{0})^2+q_0+{o(1)\/n},\ \ 
\qq \f_{0}={1\/2}\arccos {1+8(c_+-1+c_+^2)\/9}\in [0,{\pi\/2}]. 
\]
\no  iii) Let $a\in (0,{\pi\/N})$ and $n\ge 1$ be odd. Let ${jN\/3}=k_j+\wt k_j, \wt k_j\in [0,1),j=1,2, $ for some integer $k_j\ge 0$.

If $a\in \wt A=\{\wt k_1,\wt k_2\}$, then each gap $|G_n(a)|\ge 0$ and
$|G_n(a)|\to 0$  as $n\to \iy$.

If $a\notin \wt A$, then each gap $|G_n(a)|>0$ and
$|G_n(a)|\to \iy$  as $n\to \iy$.

\no iv) There exists $\ve>0$ such that:
each gap $G_n(a)\ss G_n(a')$ and each spectral bands
$S_n(a')\ss S_n(a)$ for all $0\le a<a'\le \ve, n\ge 1$.

\no iv) Let  $a\in (0,{\pi\/N})$ for even $N$
and $a\in (0,{\pi\/N})\sm \{{\pi\/2N}\}$ for odd $N$. 
Then
\[
\lb{T7-5}
S_n(a)=\cup_{k=1}^N \s_{k,n}(a),\qqq
S_n^N(a)=\cap_{k=1}^N\s_{k,n}\ne \es, \qq all\ n\ge 1,
\]
v) Let  $a_m={\pi\/2}-{\pi m\/N}\in [0,{\pi\/N}]$ for some $m\in \ol N$. Then 
\[
\lb{T7-6}
S_n(a_m)=\cup_{k=1,k\ne m}^N \s_{k,n}(a_m),\qq
\cap_{k=1,k\ne m}^N \s_{k,n}\ne \es, \qq all\ n\ge 1,
\]
and
\[
\lb{T7-7}
S_n^N(a)\to \wt \s_n=\{\wt \l_n\} \ \ \ as \qq a\to a_m.
\]
\end{theorem}

We present the plan of the paper.
In Sect. 2 we construct the fundamental solutions $\vT,\F$
and describe the basic properties of the monodromy operator.
In Sect. 3 we prove Theorem 1.2 about the eigenfunctions.
In Sect. 4 we prove the basic properties of the 
Lyapunov functions $F_k$ for the case $a=0$.
These functions are important to study the functions $F_k, k\in \Z_N, a\in \R$. In Sect. 5 we prove the basic properties of the 
Lyapunov function $F_k, k\in\ol N$. In Sect.6 we shortly recall
the results about the direct integral.

\section{ Fundamental solutions}
\setcounter{equation}{0}

In this Section we study the fundamental solutions
for the operator $H_k, k\in \Z_N$.

\begin{lemma}
\lb{mafi}
The identities \er{aoI} hold and $R_N={\sqrt3\/4\sin {\pi\/2N}}
$ is the  radius of the tube $\G^{(N)}$.
\end{lemma}
\no{\bf Proof.} Identity $\mA={B\/2}[e_0,r_\o^0+te_\o], e_0=(0,0,1)\in \R^3$ yields
\[
\label{ao2}
a_\o(t)={B\/2}([e_0,r_\o^0+te_\o],e_\o)=
{B\/2}([e_0,r_\o^0],e_\o)=a_\o(0)=a_\o,\qqq any \qq t\in[0,1].
\]
The definition $e_{n,0,k}=e_0$ gives $a_{n,0,k}=0$. Consider other cases.
Let $\o=(0,1,1)$.
Using the simple properties of the tube $\G^{(N)}$ (see Fig. 1)
we deduce that the vector $e_\o$ has the form
$e_\o=(e_\o',e_\o'', {1\/2})\in\R^3$. We will determine
 $e_\o',e_\o''$. 

The vectors
$e_\a, \a=(0,1,k), k\in \Z_N$ joint the points $r_\a^0, r_\a^1$.
The points $r_\a^0, \a=(0,1,k), k\in \Z_N$ belong to the vertices of the
N-rectangle. 
The points $r_\a^1, \a=(0,1,k), k\in \Z_N$ belong to vertices of the
another N-rectangle, turned by the angle ${\pi\/2N}$. 
Then there exists a coordinate systen in $\R^3$
such that   the vector 
$$
e_\a=(R_N\cos \f,R_N\sin \f, {1\/2})-(R_N,0,0)\in\R^3,\ \ \f={\pi\/N}.
$$  
Thus the identity $|e_\a|=1$ yields
$$
1=|e_\a|^2=R_N^2(1-\cos \f)^2+R_N^2\sin^2 \f+{1\/4},
$$
which yields $R_N={\sqrt 3\/4\sin {\f \/2}}$.
We rewrite the vector $r_\a^0$ in the form $r_\a^0=Ce_0+r_1$
for some $C\in \R$, where $r_1\bot e_0$
and $r_1=R_N(1,0,0)$. Then using \er{ao2} we obtain
$$
a_\a=
{B\/2}([e_0,r_\a^0],e_\a)={B\/2}([e_0,r_1],e_\a)={B\/2}\det
\ma R_N(\cos\f-1) & R_N\sin\f & {1\/2}\\
0 & 0 & 1\\
1 & 0 & 0\am={BR_N\/2}\sin\f,
$$
which yields $a_\a={B\sqrt3\/4}\cos {\pi\/2N}$. The proof
of other cases are similar.
\BBox

Let $H_{0k}$ be the operator $H_k$ at $q=0$.
Recall that $H_{0k}, \gD(H_{0k})=W^2(\G^{(1)})$ is  self-adjoint \cite{Ca1}.
Let $\|f\|^2=\int_{\G^{(1)}}|f(x)|^2dx$ and let $\|f_\o\|_0^2=
\int_0^1|f_\o(t)|^2dt$ for $f_\o\in L^2(\G_\o),\o\in \cZ$.
Repeating  the standard arguments from \cite{KL} for the case $a=0$ we obtain 
\[
\lb{selfa}
\|qf\|^2\le \ve \|H_{0k}f\|^2+\rt({\|q\|_0^2\/\ve}+2\|q\|_0 \rt)\|f\|^2,
\ \ all\ \ f\in \gD(H_{0k}),
\]
for any constant $0<\ve <1$. Then by the Kato-Rellich Theorem (see 162 p. \cite{ReS1}), $H_k=H_{0k}+q$ is self-adjoint on $\gD(H_{0k})$
and essentially self-adjoint on  any core of $H_{0k}$.

\no {\bf Proof of Theorem \ref{T1}.}
i) We need a modification.
We define the unitary operator $\mU$ and a modified operator
$H$ by
\[
(\mU f)_\o=e^{ita_\o}f_\o, \o\in \cZ,\ \qqq
H=\mU^* \mH \mU \ \ 
\]
in $L^2(\G^{(N)})$.
Then  acting on the edge $\G_\o$, $H$ is the ordinary differential operator given by
\[
\lb{se22}
(Hf)_\o=-f_\o''+qf_\o,\quad \ \ \ \ \ \
\ \ \  f_\o, f_\o'' \in L^2(\G_\o),\ \ \ \o\in \cZ,\qq 
\]
where  $f\in \gD(H)$ (see below) satisfies the Sturm-Liouville type of boundary  conditions

\no {\bf The Modified Kirchhoff  Boundary Conditions.} {\it Each $f\in \gD(H)$ is continuous on $\G^{(N)}\sm V^{(N)}$ and satisfies}
\[
\lb{KirC0}
f_{\o_3}(1)=f_{\o}(0)=e^{ia_\o}f_{\o_4}(1),\qqq 
f_{\o_1}(0)=e^{ia_\o}f_{\o}(1)=f_{\o_2}(0),
\]
\[
\lb{KirC1}
-f_{\o_3}'(1)+f_{\o}'(0)-e^{ia_\o}f_{\o_4}'(1)=0,\qqq
f_{\o_1}'(0)-e^{ia_\o}f_{\o}'(1)+f_{\o_2}'(0)=0,
\]
$$
all \ 
\o_1=(n+1,0,k),\ \  \o=(n,1,k),\ \o_2=(n,2,k),\
\o_3=(n,0,k),\ \ \o_4=(n,2,k-1)\in \cZ.
$$

We will obtain another representation of $H$.
The function $f$ in \er{se22}-\er{KirC1} is a vector function  $f=(f_{\o}), \o=(n,j,k)\in\cZ$.  We define a vector-valued function $f_{n,j}=(f_{n,j,k})_{k}^{N}$, where each $f_{n,j},(n,j)\in \Z\ts \Z_3$ is an $\C^N-$ vector, which satisfies the equation
\[
\lb{eqP}
-f_{n,j}''+qf_{n,j}=\l f_{n,j},
\]
and the Kirchhoff boundary conditions (which follow from \er{KirC0}-\er{KirC1})
\[
\label{C0}
f_{n,0}(1)=f_{n,1}(0)=e^{i a}\cS f_{n,2}(1),\quad
f_{n+1,0}(0)=e^{i a}f_{n,1}(1)=f_{n,2}(0),
\]
\[
\label{C1}
-f'_{n,0}(1)+f'_{n,1}(0)-e^{i a}\cS f'_{n,2}(1)=0,\quad
f'_{n+1,0}(0)-e^{ia}f'_{n,1}(1)+f'_{n,2}(0)=0,
\]
for all $n\in\Z$. Recall that the constant $a\in\R$ is given by  \er{aoI} and  the unitary operator $\cS$ in $\C^N$ is defined by
$$
\cS(h_1,h_2,\ldots,h_N)^\top=(h_N,h_1,\ldots,h_{N-1})^\top,\ \ \ 
h=(h_n)_1^N\in \C^N.
$$ 
We need more detail representation of $H$. We have
\[
\cS=\sum_1^Ns^k\cP_k,\ \  \qqq  \cP_kh=e_k(h,e_k),\qqq e_k={1\/N^{1\/2}}(1,s^{-k},s^{-2k},...,s^{-kN+k}),h\in \C^N,
\]
where $e_k$ is an eigenvector of $\cS$ and $\cS e_k=s^ke_k$
and recall $s=e^{i{2\pi \/N}}$. The operators $\cS$ and $H$ commute, then
$
H=\sum_1^N\os H_k, \qq H_k=H\cP_k.\ 
$
Using \er{C0} and \er{C1}  we deduce that the operator $H_k$ on the graph $\G^{(1)}$ acts in the Hilbert space $L^2(\G^{(1)})$. 
Then acting on the edge $\G_\o$, $H_k$ is the ordinary differential operator given by
\[
(H_kf)_{n,j}=-f_{n,j}''+q(t)f_{n,j},\qqq f_{n,j}: [0,1]\to\C
\]
on the vector functions $f=(f_{n,j}), (n,j)\in\Z\ts \Z_3$, which 
satisfy the boundary conditions
\[
\label{K0}
f_{n,0}(1)=f_{n,1}(0)=e^{i a}s^k f_{n,2}(1),\qq
f_{n+1,0}(0)=e^{i a}f_{n,1}(1)=f_{n,2}(0), \qq s=e^{i{2\pi\/N}},
\]
\[
\label{K1}
-f'_{n,0}(1)+f'_{n,1}(0)-e^{i a}s^k f'_{n,2}(1)=0,\quad
f'_{n+1,0}(0)-e^{ia}f'_{n,1}(1)+f'_{n,2}(0)=0.
\]

ii) For fixed $k\in\Z_N$ we consider the system 
\[
\label{sys}
-f''+qf=\l f,\quad \l\in\C,\;
f\in C^1(\G^{(1)}\sm V^{(1)})\text{ and }
f\text{ satisfies }\eqref{K0}\text{ and }\eqref{K1}.
\]
Suppose $f$ satisfies \er{sys} and we know $f_{0,0}(0)$ and $f_{0,0}'(0)$.  Here we will detremine the solution of this system. 

Any solution $f$ of the equation $-f''+qf=\l f$ satisfies
\[
\lb{gs1}
f(t)=\vt_tf(0)+{\vp_t\/\vp_1}(f(1)-\vt_1f(0)),\qqq t\in[0,1],
\]
\[
\label{gs2}
f'(0)={f(1)-\vt_1f(0)\/\vp_1},\quad
f'(1)={\vp_1'f(1)-f(0)\/\vp_1},
\]
where $\vt_t=\vt(1,\l), \vp_t=\vp(1,\l),...$
Substituting \er{gs2} into the first Eq. in  \er{K1} at $n=0$, we obtain
$$
-\vp_1f'_{0,0}(1)+(f_{0,1}(1)-\vt_1f_{0,1}(0))
-e^{i a}s^k(\vp_1'f_{0,2}(1)-f_{0,2}(0))=0.
$$
Using \er{K0}, we obtain
\begin{multline}
\vp_1f'_{0,0}(1)+(a^{-i a} f_{1,0}(0)-\vt_1 f_{0,0}(1))
-e^{i a}s^k(\vp_1'e^{-ia}s^{-k}f_{0,0}(1)-f_{1,0}(0))=0,
\end{multline}
which implies 
\[\label{sing}
-\vp_1f'_{0,0}(1)+u f_{1,0}(0)-2\D f_{0,0}(1)=0,\qqq
u=a^{-i a}+e^{i a}s^k=2c_ks^{{k\/2}}.
\]
If we assume that $c_k\neq 0$, then  we obtain \textbf{the first basic identity}
\[
\lb{t1}
f_{1,0}(0)=
{\vp_1\/u} f'_{0,0}(1)+{2\D\/u} f_{0,0}(1),\qqq u=2c_ks^{{k\/2}},\qq
c_k=\cos (a+{\pi k\/ N}).
\]
Substituting \er{gs2} into the second Eq. in \er{K1} at $n=0$, we get
$$
\vp_1f'_{1,0}(0)-e^{ia}(\vp_1'f_{0,1}(1)-f_{0,1}(0))
+(f_{0,2}(1)-\vt_1f_{0,2}(0))=0.
$$
Using \er{K0}, we obtain
\begin{multline}
\vp_1f'_{1,0}(0)-
e^{ia}(\vp_1'e^{-i a}f_{1,0}(0)-f_{0,0}(1))
+(e^{-ia}s^{-k}f_{0,0}(1)-\vt_1f_{1,0}(0))=0,
\end{multline}
which yields
\[\label{t2}
\vp_1f'_{1,0}(0)=2\D f_{1,0}(0)-u_2 f_{0,0}(1),\qqq u_2=e^{-i a}s^{-k}+e^{i a}=2c_ks^{-{k\/2}}.
\]
Then substituting \er{t1} into \er{t2}, we obtain 
\textbf{the second basic identity}
\[
\label{psi31}
f_{1,0}'(0)={4(\D^2-c_k^2)\/u\vp_1}f_{0,0}(1)+
{2\D\/u} f_{0,0}'(1), \qqq u=2c_ks^{{k\/2}}.
\]
Substituting consequently $f=\vT_k$ and $f=\F_k$ into
\er{t1} and into \er{psi31}, and after this $\vT_{k,\o}(0,\l), \vT_{k,\o}'(0,\l),..$ into 
$\cM_k(\l)=\ma \vT_{k,\a}(0,\l) & \F_{k,\a}(0,\l)\\
\vT_{k,\a}'(0,\l) & \F_{k,\a}'(0,\l)\am,\ \a=(1,0)$, we obtain
\er{T1-1}.

Using \er{T1-1} we obtain 
$
\det \cT_k={s^{-k}4c_k^2\/4c_k^2} =s^{-k}, \  \det \cM_k=\det \cT_kM=s^{-k}$, 
since $\det M=1$. 
Moreover, using  
$
\cR \cM\cR^{-1}=\ma \vt_1&1\\
\vp_1\vt_1'&\vp_1'\am, \quad
\ \ \cR=\ma 1 & 0\\0 & \vp_1\am,
$
we obtain 
$$
\Tr\cM_k=\Tr\cR^{-1}\cT_k \cR \cM={2s^{-{k\/2}}\/c_k}\Tr \ma 2\D & 1 \\ 4\D^2-4c_k^2 & 2\D \am
\ma \vt_1&1\\
\vp_1\vt_1'&\vp_1'\am\\
$$$$
={2s^{-{k\/2}}\/c_k}\rt(2\D \vt_1+\vp_1\vt_1'+4(\D^2-c_k^2)+2\D \vp_1'\rt)
={2s^{-{k\/2}}\/c_k}
\rt(2\D^2+{\vp_1\vt_1'\/4}-c_k^2\rt)={2s^{-{k\/2}}(F+s_k^2)\/c_k}
$$
which gives \er{T1-2}.
\BBox

\section{Eigenfunctions of $H$, Proof of Theorem \ref{T2}, \ref{T3}}
\setcounter{equation}{0}

\no {\bf Proof of Theorem \ref{T2}}.
Let $\vp_1=\vp(1,\l), \vp_1'=\vp'(1,\l),...$

\no i) Here we use arguments from \cite{KL}. Let $H_k\p=\l \p$ for some eigenfunction $\p$ and some $\l\in\s_D$.
For each  $\p_\o, \o=(n,j)\in\Z\ts\Z_3$,
we have $\p_\o(t)=\p_\o(0)+\int_0^t \p_\o'(s)ds, t\in[0,1]$. Then
$|\p_\o(0)|\leq |\p_\o(x)|+\int_0^1|\p_\o'(t)|dt$ and integrating over $x\in[0,1]$, using the H\"older inequality, we obtain 
\[\label{Fto0}
|\p_\o(0)|\le \int_0^1 (|\p_\o(t)|+|\p_\o'(t)|)dt\le \|\p_\o\|_0+\|\p_\o'\|_0=o(1), \qq as \qq n\to\pm\iy,
\]
since $\p\in\cH_k(\l)$. 
Furthermore, each restriction $\p_\o$ has the form
$\p_\o(t)=a_\o\vt(t,\l)+b_\o\vp(t,\l), t\in [0,1]$
for some constants $a_\o, b_\o$,
which implies  $\p_\o(1)=\p_\o(0)\vt_1$.
The Kirchhoff conditions give
$\p_{n+1,0}(0)=\p_{n,1}(1)=\vt_1\p_{n,1}(0)=\vt_1\p_{n,0}(1)=\vt_1^2\p_{n,0}(0)$, which yields 
\[\label{Fninf}
\p_{0,0}=\vt_1^{2n}\p_{-n,0}(0),\ \ \ 
all \ n\in\Z.
\] 
Let $|\vt_1|\le 1$, the proof for the case
$|\vt_1|> 1$ is similar. Then \er{Fninf}, \er{Fto0} imply $\p_{0,0}=\vt_1^{2n}\p_{-n,0}(0)=o(1)$ as  $n\to+\iy$.
Thus \er{Fninf} gives $\p_{n,0}(1)=0$ for all $n\in\Z$.
Finally, $\p$ vanishes on all vertexes of $\G^{(1)}$, since the set of all ends of vertexes $\G_{n,0}, n\in\Z$
coincides with the vertex set of $\G^{(1)}$.

ii) Let $A=1/\e, \e=1-e^{i2a}s^k{\vp_1'}^2$.
We will construct the eigenfunction $\p^0$.
Define a function $\f$ by:
 $\f_{n,j}=0$ for all $n\ne 0,-1$ and any $j\in \Z_3$ and
$
\f_{0,0}=\vp_t,\qq  \f_{0,j}=C_1\vp_t,\ \f_{0,j}=C_2\vp_t, t\in [0,1].
$
Then the  Kirchhoff  Conditions \er{K1} at $n=0$ give
$$
-\vp_1'+C_1-e^{ia}s^k C_2\vp_1'=0,\qqq -e^{ia}C_1\vp_1'+C_2=0,
$$
which yields
$ C_1=A\vp_1',\  C_2=Ae^{ia}{\vp_1'}^2$.

Let 
$\f_{-1,0}=0,\
\f_{-1,1}=C_3\vp_t,\   \f_{-1,2}=C_4\vp_t, \ t\in [0,1].
$
Then the  Kirchhoff  Conditions \er{K1}  at $n=-1$ give
$$
1-C_3e^{ia}\vp_1'+C_4=0,\qqq
C_3-e^{ia}s^k C_4\vp_1'=0,
$$
which yields
$
C_4=-A,\ C_3=-Ae^{ia}s^k{\vp_1'}$.

 Let $\|f\|^2=\int_{\G^{(1)}}f^2(x)dx$.
Using \er{T2-1}-\er{T2-2}, we deduce that $\p^{(0)}$ is an eigenfunction of $H_k$ with the norm given by 
$$
\int_\G|\p^{(0)}(x)|^2dx=C\int_0^1|\vp(x,\l)|^2dx,\ \ C=
|\e|^2+1+{\vp_1'}^2+2{\vp_1'}^4.
$$

The operator $H_k$ is 1-periodic, then each $\p^{(n)}, n\in \Z$ is an eigenfunction. We will show that the sequence $\p^{(n)}, n\in\Z$ forms
a basis for $\cH_k(\l)$.  The functions $\p^{(n)}$ are linearly independent, since  $\G_{n,0}\ss\supp\p^{(n)}\sm \supp\p^{(m)}$
 for every $n\ne m$.

For any $f\in \cH_k(\l)$ we will show the identity \er{T2-3}, i.e.,
\[
\lb{eit2-2}
f=\wh f,\ \ \wh f=\sum_{n\in \Z}\wh f_n\p^{(n)},\qqq
\wh f_n={f_{n,0}'(0)\/\e}.
\]
From \er{eit2-2} and $\l\in \s_D$ we deduce that
\[
\lb{idFf}
\wh f|_{\G_{n,0}}=f|_{\G_{n,0}}=\wh f_n\vp
 \ \ \ all \ \ n\in\Z.
\]
We will prove the following simple properties of $\wh f$:
\[ 
\lb{relf}
 \ \wh f\in W_{loc}^2(\G^{(1)}),\qqq
 \ \sum |\wh f_n|^2<\iy, \qqq
\ \wh f\in L^2(\G^{(1)}).
\]
Each function $\p^{(n)}\in W_{loc}^2(\G^{(1)}), n\in\Z$, then we get $\wh f\in W_{loc}^2(\G^{(1)})$. 

Let $\G^0=\cup_{n\in \Z}\G_{n,0}$.
 Using the identity \er{idFf}, we obtain
$$
\int_{\G^0}|f(x)|^2dx =\int_{\G^0}|\wh f(x)|^2dx =\sum_{n\in \Z} |\wh f_n|^2\int_{\G_{n,0}}|\p^{(n)}(x)|^2dx=\|\p^{(0)}\|^2\sum_{n\in \Z} |\wh f_n|^2
$$
which yields $\{\wh f_n\}\in \ell^2$ and 
$
\sum_{n\in \Z} |\wh f_n|^2\le \|f\|^2
$.
We will show $\|\wh f\|^2=\int_{\G^{(1)}}|f(x)|^2dx<\iy$.
Using $\p^{(n)}\p^{(m)}=0$ for $|n-s|>1$ we have
$$
\bar f f=\sum_{n,m\in Y}\ol {\wh f_n} \wh f_m\ol{\p^{(n)}}\p^{(m)},\ \ 
Y=\{n,m\in \Z:\ |n-s|\le 1\},
$$
\[
\lb{idf21}
\|f\|^2=\int_{\G}\bar f(x)f(x)dx=\sum_{n,n\in Y} \ol{\wh f_n} \wh f_m(\p^{(n)},\p^{(m)})\le \sum_{n,n\in Y} |\wh f_n| |\wh f_m|\|\p^{(n)}\|\|\p^{(m)})\|
\]
$$
=\|\p^{(0)}\|^2\sum_{n,m\in Y} |\wh f_n| |\wh f_m|
\le 3\|\p^{(0)}\|^2 \sum_{n\in \Z} |\wh f_n|^2
$$
which yields $\|\wh f\|^2\le 3\|\p^{(0)}\|^2\sum_{n\in \Z} |\wh f_n|^2$.
We have proved the properties \er{relf}.

Consider the function $u=f-\wh f$. The properties \er{relf} and \er{eit2-2}
give  
$$
u\in W_{loc}^2(\G),\ \ \ u\in L^2(\G),\ \ \ 
u|_{\G^0}=0.
$$
Define the function $u_0=u|_{\G_{0,0}}$ and
$u_0|_{\G\sm\G_{0,0}}=0$. We deduce that the function $u_0$ is an eigenfunction and has a compact support $\G_{0,0}$. But the Kirchhoff boundary conditions \er{1K0}-\er{1K1} yields $u_0=0$.

iii) Let $e^{i2a}s^k{\vp_1'}^2=1$. Then  $e^{i2a}s^k=1={\vp_1'}^2$.
Let 
$
\f_{n,0}=0, \f_{0,1}=\vp_t,\qq  \f_{0,2}=C\vp_t.
$
Then the  Kirchhoff conditions \er{1K1}  at $n=0$ yield:\ 
$
1-e^{ia}s^k C\vp_1'=0,\ -e^{ia}\vp_1'+C=0,
$
and thus $C=e^{ia}\vp_1'$

The operator $H_k$ is periodic, then each $\p^{(n)}, n\in \Z$ is an eigenfunction. We will show that the sequence $\p^{(n)}, n\in\Z$ forms
a basis for $\cH_k(\l)$.  The functions $\p^{(n)}$ are linearly independent, since $\supp\p^{(n)}\cap\supp\p^{(m)}=\es$  for all $n\ne m$.

For any $f\in \cH_k(\l)$ we will show the identity \er{T2-2}, i.e.,
\[
\lb{iii1}
f=\wh f,\ \ \wh f=\sum_{n\in \Z}\wh f_n\p^{(n)},\qqq
\wh f_n=f_{n,1}'(0).
\]
From \er{iii1} and $\l\in \s_D$ we deduce that
\[
\lb{iii2}
\wh f|_{\G_{n,1}}=f|_{\G_{n,1}}=\wh f_n\vp \qq all \ \ n\in\Z.
\]
Using the arguments from ii) we obtain the following properties of $\wh f$:
\[ 
\lb{iii3}
 \wh f\in W_{loc}^2(\G^{(1)}),\qqq
 \ \sum |\wh f_n|^2<\iy, \qqq  \wh f\in L^2(\G^{(1)}).
\]

Consider the function $u=f-\wh f$. The properties \er{iii3} and \er{iii1} give  
\[
u\in W_{loc}^2(\G^{(1)}),\qqq u\in L^2(\G^{(1)}),\qqq
u|_{\G_{n,1}}=0 \ \ all \ n\in\Z.
\]
The function $u=0$ at all vertaces of $\G^{(1)}$
and then $u_\a=C_\a\vp_t, \a=(n,j), n\in \Z, j=0,2$.
Assume that $C_{0,0}=C$. Then the Kirchhoff boundary conditions
\er{1K1} yields $C_{n,0}=-C_{n,2}=-C$ and $C_{n,0}=C_{n+1,0}=C$,
which give $C=0$ since $u\in L^2(\G^{(1)})$.
\BBox

\no {\bf Proof of Theorem \ref{T3}}.
If $c_k=0$, then by \er{TA-1}, $\s(H_k(a))=\s_{\iy}(H_k(a))$.
The case $\l\in \s_D$ is described by Theorem \ref{T2}.

Consider eigenfunctions  $-\f''+q\f=\l \f$ on $\G^{(1)}$
under the condition $c_k=0, \l\notin \s_D$. Then \er{sing} and \er{t2} give 
\[
\lb{se1}
\f'_{n,0}(0)={2\D\/\vp_1} \f_{n,0}(0),\ \qqq \f'_{n,0}(1)=-{2\D\/\vp_1} \f_{n,0}(1),\ \qqq all \ n\in \Z,
\]
which under the condition $\f_{0,0}(1)=1$ yields 
\[
\lb{se2}
\f_{0,0}(t)=\vt_t C_1+\vp_t C_2,\qqq C_1=\vp_1'+2\D=f_{0,0}(0),\qqq C_2=-\vt_1'-{2\vt_1\D\/\vp_1}=f_{0,0}'(0).
\]
Then the identities \er{se2} imply
\[
\lb{se3}
{\f_{0,0}'(0)\/\f_{0,0}(0)}=-{\vt_1'+{2\vt_1\D\/\vp_1}\/\vp_1'+2\D}={2\D\/\vp_1},
\]
and then
\[
\lb{se4}
4\D^2+2\D(\vp_1'+\vt_1)+\vp_1\vt_1'=0,\qqq 2\D^2+{\vp_1\vt_1'\/4}=0.
\]
Then only for $\l\in \s_{AP}$ there exists a
 solution of $-\f''+q\f=\l \f$.

Consider an eigenfunction $-\f''+q\f=\l \f, \l\in \s_{AP}$.
Assume that $\f_{n,0}=0$ for all $n\ne 0$. Then we get
\[
\lb{se5}
1=\f_{0,0}(1)=\f_{0,1}(0)=-e^{-ia}\f_{0,2}(1),\qqq
\f_{1,0}(0)=\f_{0,1}(1)=\f_{0,2}(0)=0,\qqq
\]
which yields
\[
\lb{se6}
\f_{0,0}(t)=\vp_t C_1+\vt_t C_2,\qqq 
\f_{0,1}(t)=\vt_t-\vp_t{\vt_1\/\vp_1},\qqq 
\f_{0,2}(t)=-{\vp_t\/\vp_1}e^{ia},\qq t\in[0,1].
\]
The direct calculation shows that the 
 Kirchhoff  Boundary Conditions \er{1K1} hold.

Consider another boundary. The  Kirchhoff  Boundary Conditions \er{1K0} at $n=-1$  give 
\[
\lb{se7}
\f_{-1,0}(1)=0=\f_{-1,1}(0)=\f_{0,2}(1),\qqq
C_1=\f_{0,0}(0)=e^{ia}\f_{-1,1}(0)=\f_{-1,2}(0),
\]
which yields
\[
\lb{se8}
\f_{-1,1}(t)={\vp_t\/\vp_1}e^{-ia}C_1,\qqq 
\f_{-1,2}(t)=(\vt_t-\vp_t{\vt_1\/\vp_1})C_1,\qqq C_1=\vp_1'+2\D.
\]
The direct calculation shows that the 
 Kirchhoff  Boundary Conditions \er{1K1} hold at $n=-1$.

 The operator $H_k$ is periodic, hence
each $\p^{(n)}, n\in \Z$ is an eigenfunction.
We will show that the sequence $\p^{(n)}, n\in\Z$ forms
a basis for $\cH_k(\l)$. The functions $\p^{(n)}$ are linearly independent, since $\supp\p^{(n)}\sm \supp\p^{(m)}=\G_{n,0}$ for every $n\ne m$.

For any $f\in \cH_k(\l)$ we will show the identity \er{T3-4}, i.e.,
\[
\lb{se9}
f=\wh f,\ \qqq \wh f=\sum_{n\in \Z}\wh f_n\f^{(n)},\qqq
\wh f_n={f_{n,0}(0)\/C_2}.
\]
The function $f$ satisfies \er{se1}, since $F(\l)=-1$.
From \er{se1} and \er{se2} we deduce that
\[
\lb{se10}
f_{0,0}(t)=C(\vt_t C_1+\vp_t C_2)=C\f_{0,0}(t),\qqq C_1=\vp_1'+2\D,\qqq C_2=-\vt_1'-{2\vt_1\D\/\vp_1},
\]
which yields
\[
\lb{se11}
\wh f|_{\G_{n,0}}=f|_{\G_{n,0}}=\wh f_n\f^{(n)}
 \ \ \ all \ \ n\in\Z.
\]
Using the arguments from the proof of  Theorem \ref{T2} we obtain the following  properties of $\wh f$:
\[ 
\lb{se12}
 \ \wh f\in W_{loc}^2(\G^{(1)}),\qqq
 \ \sum |\wh f_n|^2<\iy, \qqq
\ \wh f\in L^2(\G^{(1)}).
\]
Consider the function $u=f-\wh f$. The properties \er{se12} and \er{se11}
give  
\[
u\in W_{loc}^2(\G^{(1)}),\ \ \ u\in L^2(\G^{(1)}),\ \ \ 
u|_{\G^0}=0.
\]
Define the function $u_0=u|_{\G_{0,0}}$ and
$u_0|_{\G\sm\G_{0,0}}=0$.
We deduce that the function $u_0$ is an eigenfunction
and has a compact support $\G_{0,0}$. But the Kirchhoff boundary conditions yields $u_0=0$.
\BBox

\section{The magnetic fields $B=0$}
\setcounter{equation}{0}

In this section we consider only  the case $a=0$. Each integer $N$ has the form 
$N=2m+1$ or $N=2m+2, m\ge 0$. Recall that the results for the case $a=0$ and $N$ odd was obtained in \cite{KL}.

Let $\t_{k,\pm}$ be the eigenvalues of $\cM_k, k\in \Z_N, k\ne {N\/2}$. Using \er{T1-2} we deduce that
\[
\lb{41}
\t_{k,-}\t_{k,+}=s^{-k}, \qq \t_{k,-}+\t_{k,+}=\Tr 
\cM_{k}=s^{-k}\Tr\cM_{-k},
\qq \t_{-k,+}=\t_{k,-}s^{k},\qq \t_{-k,-}=\t_{k,+}s^{k}.
\]
If $k=0$ (or $k=N$), then we can
introduce the standard  entire Lyapunov function 
$$
\D_0=F={\Tr \cM_0\/2}={1\/2}(\t_{0,+}+{1\/\t_{0,+}})=\D^2+{\vp_1\vt_1'\/4}-1.
$$
If $k\neq 0$, then we define the Lyapunov functions 
$\D_{k,\pm}={1\/2}(\t_{k,\pm}+{1\/\t_{k,\pm}})$ (see \cite{BBK},\cite{CK}) and
using  \er{41} we get $\D_{-k,\pm}=\D_{k,\pm}$.
Below we prove the following identities
\[
\lb{0DeL1}
\D_{k,\pm}={1\/2}(\t_{k,\pm}+{1\/\t_{k,\pm}})=\x_k\pm\sqrt{\rho_k},\qqq \x_k=\D_0+s_k^2,\ \ \ \r_k={s_{0k}^2\/c_{0k}^2}(c_{0k}^2-\x_k^2),\qq k\ne {N\/2}.
\]
 If $q=0$, then we denote the corresponding functions by
$\D_k^0, \r_k^0,..$. In particular, we have
\[
\D_0^0(\l)={9\cos 2\sqrt{\l}-1\/8},\ \qqq  \x_k^0=\D_0^0+s_{0k}^2,\ \ \ \r_k^0={s_{0k}^2\/c_{0k}^2}(c_{0k}^2-(\x_k^0)^2).
\]
Introduce the two sheeted Riemann surface $\gR_k$
(of infinite genus) defined by $\sqrt {\r_k}$.
The functions $\D_{k,\pm},k\in \ol m,$ are the branches of $\D_k=\x_k+\sqrt{\r_k}$ on the Riemann surface $\cR_k$.

\begin{theorem}
\lb{T03} 
If $a=0, k={N\/2}\in \Z$, then following identities hold:
\[
\lb{T07-1}
\s(H_k)=\s_{\iy}(H_k)=\s_D\cup \s_{AP}.
\] 
Moreover, if $k=0,1,...,m,k\ne {N\/2}$, then 

\no i) The Lyapunov functions $\D_{k,\pm}, k\in\ol m$ satisfy \er{0DeL1}. 
 
\no ii) For each $k=0,1,...,m$ the  following identities hold:
\[
\lb{T03-4}
\s(H_k)=\s_{\iy}(H_k)\cup\s_{ac}(H_k),\ \ \ \ 
\s_{\iy}(H_k)=\s_D,\ \ \ \ \s_{ac}(H_k)=\{\l\in\R: \D_k(\l)\in [-1,1]\}.
\] 
iii) If $\D_k(\l)\in (-1,1)$ for some $\l\in \R,k=0,..,m$ and 
$\l$ is not a branch point of $\D_k$, then $\D_k'(\l)\neq 0$. 

\end{theorem}

{\bf Remark.} 1) If we know $\D_0$, then we determine all
$\D_k, \r_k,k=-m,..,m$ by \er{0DeL1}. 2) If we know $\r_j$ for some $j=1,..,m$, then we determine all
$\D_k, \r_k,k=-m,..,m$ by \er{0DeL1}. 

\no {\bf Proof.} 
 The proof of $k={N\/2}\in\Z$ follows from Theorem \ref{T3}. 
 The case $k=0,1,...,m$ and $N$ is odd was proved in \cite{KL}.
The proof for even $N$ is similar. \BBox

Recall that a zero of $\r_k=R_k|_{a=0}, k\in \ol m$ is called a {\bf resonance} of $H_k$.

Let $D_0(\cdot,\t)=\det (\cM_0-\t I_2),\t\in \C$ and let $\l_{0,n}^\pm$ be the zeros of $D_0(\l,1)D_0(\l,-1)$.
The zeros of $D_0(\l,1)$ ( and $D_0(\l,-1)$ counted
with multiplicity are the periodic (anti-periodic) eigenvalues
for the equation $-y''+qy=\l y$ on $\G$  with periodic (anti-periodic) boundary conditions.
Below we will show that $\l_{0,n}^\pm$ satisfy 
\[
\lb{epa0}
\D_0(\l_{0,n}^\pm)=(-1)^n,\ \ \ 
\l_{0,0}^+<\l_{0,1}^-< \l_{0,1}^+<\l_{0,2}^-\le \l_{0,2}^+<\l_{0,3}^-< \l_{0,3}^+<\l_{0,4}^-\le \l_{0,4}^+<....
\]
Introduce the gaps $\g_{0,n}=(\l_{0,n}^-,\l_{0,n}^+), n\ge 1$. 
We need the following results from \cite{KL} about $F=\D_0$.

\begin{theorem}
\label{T04}
Let $a=0$. Then the function $F=\D_0$ is entire and has the following properties:

\no i) The function $\D_0'$ has only real simple zeros
$\l_{0,n},n\ge 1$,
which are separated by the simple zeros $\e_{0,n}$ of $\D_0$:
$\e_{0,1}<\l_{0,1}<\e_{0,2}<\l_{0,2}<\e_{0,3}<...$ and satisfy
\[
\lb{T04-1}
\n_n,\m_n\in \g_{0,2n},\ \ 
\D_0(\l_{0,2n})\ge 1,\ \ \ \   \  \ \ \ \ \ \ \   \D_0(\l_{0,2n-1})\le -{5\/4},\ \ \ \ \  for \ any \ n\ge 1.
\]
\no ii) The periodic and anti-periodic eigenvalues $\l_{0,n}^\pm,n\ge 0$ satisfy \er{epa0} and have asymptotics
\[
\lb{T04-2}
\l_{0,2n+1}^\pm=(\pi{2n+1\/2}\pm \f)^2+q_0+{o(1)\/n},\ \ \qq
\  \l_{2n}^\pm=(\pi n)^2+q_0\pm \rt|{|\hat q_n|^2}-{\hat q_{sn}^2\/9}\rt|^{1\/2}+{O(1)\/n},
\]
 $\f=\arcsin {1\/3}\in [0,{\pi\/2}]$ and $\hat q_n=\int_0^1q(t)e^{i2\pi nt}dt, \hat q_{sn}=\Im \hat q_{n}$. Moreover, 
 $(-1)^n\D_0(\l)\ge 1$ for all $\l\in [\l^-_{0,n},\l^+_{0,n}]$ and 
 $\l_{0,n}\in [\l^-_{0,n},\l^+_{0,n}]$.
 
 \no iii) $|\g_n|=|\m_n-\n_n|$  iff $\g_n=\g_{0,2n}$.

\no iv) $\g_{0,2n}\ss \g_n$ for all $n\ge 1$. Moreover,
for fixed $n\ge 1$ we have $\g_{0,2n}=\es$ iff $\g_{n}=\es$.

\no v) $\D_0(\l_{1+2n})=-{5\/4}$ for all $n\ge 0$ iff 
$q\in L_{even}^2(0,1)$.

\end{theorem}

{\bf Consider the case} $k\in \ol m, m\ge 1$.
Let $D_k(\cdot,\t)=\det (\cM_k-\t I_2),\t\in \C$.
Let $\l_{k,2n}^\pm$ and $\l_{k,2n+1}^\pm,n\ge 0$ be the zeros of $D_k(\l,1)$ and $D_k(\l,-1)$. Below we will show that the periodic  eigenvalues $\l_{k,2n}^\pm$ and the anti-periodic eigenvalues $\l_{k,2n+1}^\pm$ satisfy the equations $\D_k(\l_{k,n}^\pm)=(-1)^n$.
In Theorem \ref{T05} we show that $\l_{k,n}^\pm$ satisfy the equations
\[
\lb{0eqpk}
\D_0(\l_{k,2n}^\pm)=\cos {2\pi k\/N},\ \qq \D_0(\l_{k,2n+1}^\pm)=
-1,\ \ \ \ \ \ k\in \ol m,
\]
and labeling is given by: each $\l_{k,n}^\pm$ is double and
\[
\lb{0espk}
\l_{0,0}^+<\l_{1,0}^+<\l_{2,0}^+<...<\l_{m,0}^+<\l_{0,1}^-
<\l_{0,1}^+
<\l_{m,2}^-<\l_{m-1,2}^-<...<\l_{0,2}^-\le \l_{0,2}^+<...,
\]
\[
\lb{0espkI}
\l_{0,2n-1}^\pm=\l_{k,2n-1}^\pm,\ \ (k,n)\in \ol  m\ts\N. 
\]
The periodic eigenvalues $\l_{k,2n}^{\pm}, k=1,..,m, n\ge 0$ (i.e., $\D_k(\l_{k,2n}^{\pm})=1$ ) satisfy
\[
\lb{0Tas-2}
\l_{k,2n}^{\pm}=(\pi n\pm \f_k)^2+q_0+{o(1)\/n},\ \ \qq 
\f_k={1\/2}\arccos {1+8c_{2k}\/9}\in [0,{\pi\/2}] \qq
as \ \ \ \ n\to \iy.
\]

Let $r_{k,n}^\pm, k\in \ol m, n\ge 0$ be the zeros of $\r_k$.
Note that there are no resonances for the cases $N=1,2$.
Below we will show that these resonances satisfy the equations
\[
\lb{0eqr}
\D_0(r_{k,2n}^\pm )=c_k-s_k^2, \  \qq 
\D_0(r_{k,2n+1}^\pm )=-c_k-s_k^2, \ \qq
\ \ \ k\in \ol m,
\]
they are real and labeling is given by
\[
\lb{0esr}
r_{k,0}^+<
r_{k,1}^-<r_{k,1}^+< r_{k,2}^-<r_{k,2}^+ <..,\ \ \ \ k\neq {N\/3},
\]\[
\lb{0esrI}
r_{k,0}^+<
r_{k,1}^-\le r_{k,1}^+< r_{k,2}^-<r_{k,2}^+ <r_{k,3}^-\le r_{k,3}^+..,\ \ \ \ k= {N\/3}.
\]
The resonances $r_{k,n}^{\pm}$ (for $k\in \ol m,k\ne {N\/3}$) satisfy 
\[
\lb{0Tas-4}
r_{k,n}^\pm=({\pi n\/2}\pm b_{k,n})^2+q_0+o(n^{-1}),\qq 
 as \qq n\to \iy,\qq
 b_{k,2n}=\f_{k,0},\qq b_{k,2n+1}={\pi\/2}-\f_{k,1}, 
\]
$\f_{k,s}={1\/2}\arccos {1+(-1)^s8c_{k}-8s_k^2\/9}\in [0,{\pi\/2}], s=0,1$. We describe the spectral properties of $H_k$.

\begin{theorem}\lb{T05} 
Let $a=0$.

\no i) The periodic and anti-periodic eigenvalues $\l_{k,n}^\pm, 
k\in \ol m,n\ge 0$ satisfy Eq. \er{0eqpk} and the relations \er{0espk}-
\er{0Tas-2}.

\no ii) The resonances  $r_{k,n}^\pm,k\in \ol m,n\ge 0$ satisfy Eq. \er{0eqr} and estimates \er{0esr}-\er{0Tas-4}.

\no iii) For $ n\ge 1, k\in \ol m$ the  following identities are fulfilled:
\[
\lb{T05-1}
\s_{ac}(H)=\cup_{n\ge 1} S_n=\cup_0^N \s_{k},
\ \  S_n=[E_{n-1}^+,E_n^-]=\cup_{k=0}^m \s_{k,n},\ \  \ 
\cap_{k=0}^m \s_{k,n}\ne \es,
\]
\[
\lb{T05-2}
\s_{0,n}=[\l_{0,n-1}^+,\l_{0,n}^-],\ \ \  {\rm and }\ \ \ 
\l_{k,n}^\pm\in \s_{k,n}=[r_{k,n-1}^+r_{k,n}^-],\ \ \ 
\]
\[
\lb{T05-3}
 G_n=(E_{n}^-,E_n^+)=\cap_0^m \g_{k,n},\ \ G_{2n}=\g_{0,2n}.
\]
\no iv) $G_{2n}=\es$ iff $\g_{n}=\es$. Moreover, 
$G_{2n}=o(1)$ as $n\to \iy$.

\no v) If $p={N\/3}\in \Z$, then $
G_{2n+1}=\es$ for all $n\ge 0$ iff $q\in L_{even}(0,1)$. Moreover, 
 each odd gap $G_n$ has the form $G_n=(r_{p,n}^-,r_{p,n}^+)$ ($n\ge 1$ is odd) and $r_{p,n}^\pm$ satisfy  
\[
\lb{T05-4}
r_{p,n}^\pm=\pi^2{n^2\/4}+q_0\pm |\wt q_{cn}|+o(n^{-1})\ \  as\qq n\to \iy,
\qq \wt q_{cn}=\int_0^1q(t)\cos \pi nt dt.
\]
\no vi) If ${N\/3}\notin \Z$, then each gap
$G_{2n+1}\ne \es, n\ge 0$ and $|G_{2n+1}|\to \iy$ as $n\to \iy$.

\no vii) The operator $H$ has only a finite number of non degenerate gaps $G_n$ iff ${N\/3}\in \Z$ and $q$ is an even
finite gap potential for the operator $-y''+qy$ on the real line.

\end{theorem}
\no {\bf Proof.} The case $N$ odd was proved in \cite{KL}.
The proof for the even $N$ is similar. \BBox

Note that in the matrix case (see \cite{BBK}) 
the resonances in general have non-zero imaginary part.
That all resonances are  real is a peculiarity of the high
symmetry of a zigzag periodic graph.

Below we need following results from \cite{KL}.

\begin{lemma}
\label{TLas} Let $c\in [-{5\/4},1]$. Then Eq. $F(\l)=c,\l\in \C$
has only real zeros, which satisfy

\no i) If $c\in (-{5\/4},1)$. Then these zeros $z_n^\pm$ are given by
\[
\lb{TLas-1}
z_{0}^{+}<z_{1}^{-}<z_{1}^{+}<z_{2}^{-}<.. \qq and \qq
\sqrt{z_{n}^{\pm}}=u_{n}^{\pm}+{q_0\/2u_n^{\pm}}+{o(1)\/n^2},\ \ as \ \ \ \ n\to \iy.
\]
where $u_n^{\pm}=\pi n\pm u_0^{+}, \ n\ge 1$ and $u_{0}^{+}={\arccos {1+8c\/9}\/2}\in [0,{\pi\/2}],u_{0}^{+}<u_{1}^{-}
<u_{1}^{+}<u_{2}^{-}<.. $.

Moreover, the zeros $z_n^\pm$ have another representation given by: $x_n^-=z_{n+1}^{+}, x_n^+=z_{n}^{-},$  and
\[
\lb{TLas-2}
x_{1}^{-}<x_{1}^{+}<x_{2}^{-}<x_{2}^{+}<x_{3}^{-}<.. \qqq
\sqrt{x_{n}^{\pm}}=v_{n}^{\pm}+{q_0\/2v_n^{\pm}}+{o(1)\/n^2}\qq as \qq n\to \iy, 
\]
where $v_n^{\pm}=\pi (n-{1\/2})\pm ({\pi\/2}-u_0^{+})$.

\no ii) If $c=1$. Then these zeros $z_n^\pm$ are given by
\[
\lb{TLas-3}
z_{0}^{+}<z_{1}^{-}\le z_{1}^{+}<z_2^-\le z_2^+<.. \qq and \qq
\sqrt{z_{n}^{\pm}}=\pi n+{q_0\/2\pi n}+{o(1)\/n^2},\ \ as \ \ \ \ n\to \iy.
\]
\no iii) If $c=-{5\/4}$. Then these zeros $x_n^\pm$ are given by
\[
\lb{TLas-4}
x_{1}^{-}\le x_{1}^{+}<x_2^-\le x_2^+<.. \qq and \qq
\sqrt{x_{n}^{\pm}}=\pi (n-{1\/2})+{q_0\/2\pi n}+{o(1)\/n^2},\ \ as \ \ \ \ n\to \iy.
\]
\end{lemma}

\section{Lyapunov functions, Proof of Theorem \ref{T4}-\ref{T7}}
\setcounter{equation}{0}

We recall that $\vt(x,\l)$, $\vp(x,\l)$, $x\in \R$ are the fundamental solutions  of Eq. $-y''+qy=\l y,\ \l\in\C,$ on the real line
such that $\vt(0,\l)=\vp'(0,\l)=1$, $\vt'(0,\l)=\vp(0,\l)=1$. 
For each $x\in \R$ the functions $\vt,\vt',\vp,\vp'$ are entire in $\l\in\C$. Moreover, the following asymptotics are fulfilled:
$$
 \vt(x,\l)=\cos\sqrt{\l}x+O\lt({e^{|\Im\sqrt{\l}|x}\/\sqrt{\l}}\rt),\qqq
\lb{vpas} \vp(x,\l)= {\sin\sqrt{\l}x\/\sqrt{\l}}+
O\lt({e^{|\Im\sqrt{\l}|x}\/|\l|}\rt),
$$
\[
\label{Das} \D(\l)= \cos\sqrt{\l}+
{q_0\sin \sqrt{\l}\/2\sqrt{\l}}+O\lt({e^{|\Im\sqrt{\l}|}\/|\l|}\rt),\ \ 
q_0={1\/2}\int_0^1q(t)dt,
\]
as $|\l |\to\iy$, uniformly on bounded sets of $(x;q)\in [0,1]\ts L^2_C(0,1)$  (see \cite{PT}).

A great number of  papers is devoted to the inverse spectral
theory and a priori estimates for the Hill operator: \cite{M}, \cite{GT},\cite{KK}, \cite{Ko} etc. 
We recall some needed below results about the Hill operator from \cite{M},...
The sequence $\l_0^+<\l_1^-\le \l_1^+\ <.....$
is the spectrum of equation $-y''+qy$ with periodic
boundary conditions of period 2,  that is  $f(x+2)=f(x), x\in \R$.
Here equality means that $\l_n^-=\l_n^+$ is an eigenvalue of multiplicity 2. Note that $\D(\l_{n}^{\pm})=(-1)^n, \  n\ge 1$. The lowest  eigenvalue $\l_0^+$ is simple, $\D(\l_0^+)=1$, and the
corresponding eigenfunction has period 1. The eigenfunctions
corresponding to $\l_n^{\pm}$ have period 1 if $n$ is even,
and they are anti-periodic, that is $f(x+1)=-f(x),\ x\in \R$, if
$n$ is odd. The derivative of the Lyapunov function has a zero $\l_n$ in the ''closed gap'' $[\l^-_n,\l^+_n]$, that is $ \D'(\l_n)=0$. 
Recall that $\m_n$ and $\n_n$ are the Dirichlet and Neumann eigenvalues.  
It is well-known that $\m_n, \n_n \in [\l^-_n,\l^+_n ]$ and $\n_0\le 0$.
Moreover, a potential $q$ is even, i.e., $e\in L^2_{even}(0,1)$ iff $|\g_n|=|\m_n-\n_n| $ for all $n\ge 1$ (see \cite{GT} or \cite{KK} or \cite{Ko}).

Similarly to the case of periodic Schr\"odinger operators on $\R$,
we determine the spectrum of $H$ in terms of the Lyapunov function
$F_k=T_k+\sqrt{R_k}$.

\no {\bf Proof of Theorem \ref{T4}.}
i)  Fix a eigenvalue $\t_k$ of $\cM_k$. 
Due to \er{T1-1},\er{T1-2} the characteristic equation for
$\cM_k$ has the form 
$
\det(\cM_k-\t I_2)=\t^2-\Tr \cM_k\t+s^{-k}=0.
$
Then eigenvalues $\t_k^\pm$ of $\cM_k$ are given by
\[
\t_k^\pm={s^{-{k\/2}}}(w_k\pm i\sqrt{1-w_k^2}),\qqq
w_k={F+s_k^2\/c_k}={s^{{k\/2}}\/2}\Tr \cM_k.
\]
Each eigenvalue $\t_k^\pm$ defines the Lyapunov function $F_k^\pm$    by $F_k^\pm={1\/2}\rt(\t_k^\pm+\frac{1}{\t^\pm_k}\rt)$.
Due to $\t_k^+\t_k^-=s^{-k}$, we obtain
\begin{multline}
F_k^\pm={\t_k^\pm+s^k \t_k^\mp \/2}
={1\/2}\rt({w_k\pm i\sqrt{1-w_k^2}\/s^{k\/2}}+(w_k\mp i\sqrt{1-w_k^2})s^{k\/2}\rt)\\
={s^{k\/2}+s^{-{k\/2}}\/2}w_k\pm i{s^{-{k\/2}}-s^{{k\/2}}\/2}\sqrt{1-w_k^2}
=c_{0k}w_k\pm s_{0k}\sqrt{1-w_k^2},
\end{multline}
where $s_{0k}=\sin \frac{k\pi}{N}, c_{0k}=\cos \frac{k\pi}{N}$. Finally,
$$
F_k^\pm=T_k\pm\sqrt{R_k},\qq T_k={c_{0k}\/c_k}(F+s_k^2),
\qq R_k=s_{0k}^2\rt(1-w_k^2\rt)=
{s_{0k}^2\/c_k^2}\rt(c_k^2-(F+s_k^2)^2\rt).
$$

ii) The standard arguments (see \cite{Ca1} and Teorem \ref{TA}) yield
$\s_{sc}(H_k)=\es$ and $\s_{ac}(H_k)=\{\l\in \R: F_k(\l)\in [-1,1]\}$.

iii) Let $F_k'(\l_0)=0$
and $F_k(\l_0)\in (-1,1)$ for some $\l_0\in \s_k, k=0,..,N$. Then we have the Tailor series
$
F_k(\l)=F_k(\l_0)+t^p{F_k^{(p)}(\l_0)\/p!}+O(t^{p+1})$,\
as $ t=\l-\l_0\to 0$,
where $F_k^{(p)}(\l_0)\ne 0$ for some $p\ge 1$. By the
Implicit Function Theorem, there exists some curve $Y\ss
\{\l:|\l-\l_0|<\ve \}\cap \C_+, Y\neq \es$,  for some $\ve>0$
such that $F_k(\l)\in (-1,1)$ for any $\l\in Y$. Thus we have
a contradiction with \er{T3-4}.
\BBox

{\bf Proof of Theorem \ref{T5}}.
Recall the identities
\[
\lb{pT51}
F_0={1\/2}\Tr\cM_0={F+s_0^2\/c_0},\qq F=2\D^2+{\vp(1,\cdot)\vt'(1,\cdot)\/4}-1={9\D^2-\D_-^2-5\/4},
\]
where $\D_-={1\/2}(\vp'(1,\cdot)-\vt(1,\cdot))$.
Recall that we consider the case: $a\in (0,{\pi\/2})$

i) The results of i) follow from \er{pT51}  and Theorem \ref{T04}.

ii) We determine the equation for periodic eigenvalues
for the case $c_k>$ and $a\in (0,{\pi\/2})$.
We have
$$
F_0(\l)=\pm 1 \lra F(\l)=\pm c_0-s_0^2,
$$
where
$
c_0-s_0^2\in (-1,1),\ \ -c_0-s_0^2\in [-{5\/4},-1).
$
Then the properties of $F$ from Theorem \ref{T04} give that all even gaps $\g_{0,2n}(a)=(\l_{0,2n}^-(a),\l_{0,2n}^+(a))$ are open  and satisfy 
$\g_{0,2n}(a)\ss \g_{0,2n}(a_1),\ a<a_1<{\pi\/2}$.
 
Consider the odd gaps $\g_{0,m}(a)=(\l_{0,m}^-(a),\l_{0,m}^+(a)), m=2n-1$. Then the properties of $F$ from Theorem \ref{T04} give
$$
\g_{0,m}(a)\supset \g_{0,m}(a_1), \qqq 0\le a< a_1<{\pi\/3},\qq and
\qq
 \g_{0,m}(a)\ss \g_{0,m}(a_1), \qqq {\pi\/3}\le a<a_1<{\pi\/2}.
 $$
If $a={\pi\/3}$, then $c_0+s_0^2={5\/4}$ and  
the properties of $F$ from Theorem \ref{T04} give
that $|\g_{0,m}(a)|\ge 0$. Moreover, if in addition,
$q\in L_{even}^2(0,1)$, then all odd $\g_{0,m}(a)=\es, m\ge 1$.

iii) Repeating the standard arguments (see \cite{KL}) we determine  asymptotics \er{Tas-4}.

iv) Consider the case: $c_0\downarrow  0$ or $a\uparrow{\pi\/2}$.
The spectral bands are defined by
$$
\s_{0,n}(a)=(\l_{0,n-1}^+(a),\l_{0,n}^-(a)), \qqq F(\l_{0,2n}^\pm)=c_0-s_0^2,\qqq
F(\l_{0,2n-1}^\pm)=-c_0-s_0^2.
$$
Using $F(\wt\l_{n})=-1$ we obtain 
$$
F(\l_{0,n}^-(a))-F(\wt\l_{n})=c_0-s_0^2+1=c_0(1-c_0)\to 0,
$$
and $F'(\wt\l_{n})\ne 0$ yields
$
\l_{0,n}^-(a)=\wt\l_{n}+{c_0 +O(c_0 ^2)\/F'(\wt\l_{n})}
\qq as \qq a\uparrow{\pi\/2}.
$
Similar arguments give
$\l_{0,n-1}^-=\wt\l_{0,2n+1}^-+c_0(1+O(c_0 ))/F'(\wt\l_{n})$
 as $a\uparrow{\pi\/2}$. Thus we have \er{T5-4}.
\BBox

We consider the Lyapunov function $F_k, k\ne N$.

{\bf Proof of Theorem \ref{T6}}.
Recall
$$
F_k^\pm=T_k\pm\sqrt{R_k},\qqq
T_k={c_{0k}\/c_k}(F+s_k^2)
\qq R_k={s_{0k}^2\/c_k^2}\rt(c_k^2-(F+s_k^2)^2\rt),\qq
k\in \ol{N-1},
$$
where $c_{0k}=\cos {\pi k\/N}, s_{0k}=\sin {\pi k\/N}$.
For the case $c_{k}\ne 0$ these definitions yield :
\[
T_k(\l,a+\pi)=-T_k(\l,a),\qq R_k(\l,a+\pi)=R_k(\l,a),\qq
\l\in \C.
\]

ii) We determine the equation for resonances for the case $k\in \ol{N-1}$. 
If $R_k(\l)<0$, then $F_k(\l) $ is not real and then $\l\notin \s(H_k)$.
We have
\[
\lb{512}
R_k<0\qq  \lra  \qq c_k^2<(F+s_k^2)^2\qq \lra \qq \ F+s_k^2>|c_k|\qq or
\qq F+s_k^2<-|c_k|.
\]
Consider the first case: $n$ is {\bf even}, i.e., $r_{k,2n}^\pm, n\in \Z$. We have
\[
F+s_k^2>|c_k|\ \ \ \ \  \lra  \ \  F>f_k^+=|c_k|-s_k^2\in [-1,1],\qqq k\in \ol{N-1}.
\]
The resonances $r_{k,n}^\pm$  are zeros of Eq. $F(\l)=f_k^+\in [-1,1],\l\in \C$. Lemma \ref{TLas} gives that all these resonances  are real and labeling is given by  \er{esr}-\er{esrI}  and $r_{k,n}^\pm $ satisfy asmptotics \er{Tas-4}.

Consider the second case: the {\bf odd} gaps $\g_{k,n}=(r_{k,n}^-,r_{k,n}^+)$, i.e., $n\in \Z$ is odd. We have
\[
 F+s_k^2<-|c_k|\qqq  \lra  \qqq
F<f _k^-=-|c_k|-s_k^2\in [-{5\/4},-1),\qq k\in \ol{N-1}.
\]
The resonances $r_{k,n}^\pm$  are zeros of  Eq. $F(\l)=f_k^-\in [-{5\/4},-1),\l\in \C$. 
Lemma \ref{TLas} gives that all these resonances  are real and labeling is given by  \er{esr}-\er{esrI}  and $r_{k,n}^\pm $ satisfy 
asmptotics \er{Tas-4}.
In order to prove \er{T6-2} below we will show  that if $R_k(\l)\ge 0$, then $F_k(\l)\in [-1,1]$.

iii) Using \er{T1-1},\er{T1-2}, we obtain the equation for periodic (anti-periodic) eigenvalues
$$
0=\det (\cM_k \mp I_2)=1 \mp \Tr M_k+s^{-k}=
1+s^{-k}\mp 2{F+s_k^2\/s^{k\/2}c_k}=2s^{-k\/2}
\rt(c_{0k}\mp {(F+s_k^2)\/c_k}\rt),
$$
  which yields the equations  
\[
\lb{510}
F(\l_n^\pm)=c_{0k}c_k-s_k^2\in (-{5\/4},1) \ \ \ \ {\rm for \ periodic \ eigenvalues}\ \ \l_n^\pm,n\ge 0,
\]
\[
\lb{511}
F(\m_n^\pm)=-c_{0k}c_k-s_k^2 \in (-{5\/4},1) \  \ \ \ \ {\rm for \ anti-periodic \ eigenvalues}\ \ \m_n^\pm,n\ge 0.
\]
These equations and Lemma \ref{TLas} give that all
 these eigenvalues  are real and labeling is given by \er{epak1}-\er{epak3}  and they satisfy asmptotics \er{Tas-4}.
Note that $\l_n^\pm, \m_n^\pm$ are simple since $F'(\l_n^\pm)\ne 0, F'(\m_n^\pm)\ne 0$, see Theorem \ref{T04}.

We will show $\m_n^\pm,\l_n^\pm\in \s_{k,n}$ and \er{T6-2}. Recall 
\[
F_k=c_{0k}w_k+s_{0k}\sqrt{1-w_k^2},\qqq w_k={F+s_k^2\/c_k},\qqq
w_k(\l_n^\pm)=c_{0k},\qqq w_k(\m_n^\pm)=-c_{0k}.
\]
Then at $\l=\l_n^\pm$  we get
\[
F_k'(\l)=w_k'(\l)\rt(c_{0k}-s_{0k}{w_k(\l)\/\sqrt{1-w_k^2(\l)}}\rt)=0,\qq
F_k''(\l)\ne 0.
\]
and similar arguments yield $F_k'(\m_n^\pm)=0,F_k''(\m_n^\pm)\ne 0$.
This shows that there are no classical gaps, $\m_n^\pm,\l_n^\pm\in \s_{k,n}$ and there exist only resonance gaps.
 Thus we obtain \er{T6-2}.

i) Using \er{512} we obtain 
$
s(H_k(a))=\{\l\in \R: (F(\l)+s_k^2)^2\le c_k^2\}\cup \s_D,\ k\in \Z_N.
$
and 
$
\s(H_k(a+{\pi\/N}))=\s(H_{k+1}(a))\cup \s_D,\qq (k,a)\in \Z_N\ts
\R$,
which yields \er{T6-1}.

iv) The proof of the case $\s_{k,n}\to \wt\s_n$ as $a\to a_{k,m}$ repeats the case $\s_{0,n}\to \wt\s_n$ as $a\uparrow {\pi\/2}$,
see the proof of Theorem \ref{T5}.
\BBox

{\bf Proof of Theorem \ref{T7}}.
i) Using \er{T6-1} we get \er{T7-1}.

ii) Consider the case $a\in(0,{\pi\/N})$ and the gaps $G_n$ for 
even $n$. Recall that for 
$\l\in \g_{k,n}$ we have obtained in the proof of Theorem \ref{T6}
\[
F(\l)+s_k^2>|c_k|\ \ \ \ \  \lra  \ \  F(\l)>f_k^+=|c_k|-s_k^2\in [-1,1],\qqq k\in \ol{N-1}.
\]
The resonances $r_{k,n}^\pm$  are zeros of Eq. $F(\l)=f_k^+\in [-1,1],\l\in \C$. Then we deduce that
$$
1>f_0^+>f_1^+>f_2^+>...>f_{p_0}^+,\qqq 1>f_{N-1}^+>f_{N-2}^+>...>f_{p_0-1}^+,
$$ 
where ${\pi\/2}-a-{\pi p_0\/N}\in (0,1)$ for some integer $p_0\ge 0$,  which yields for $n=0$
$$
\l_{0,0}^+<r_{1,0}^+<r_{2,0}^+<...<r_{p_0,0}^+,\qqq 
r_{N-1,0}^+<r_{N-2,0}^+<...<r_{p_0-1,0}^+,
$$
 for $n\ge 2$
$$
r_{p_0,n}^-<r_{p_0-2,n}^-<...<r_{1,n}^-<\l_{0,n}^-< \l_{0,n}^+<r_{1,n}^+<r_{2,n}^+<...<r_{p_0,n}^+,
$$
$$
r_{p_0+1,n}^-<r_{p_0+2,n}^-<...<r_{N-1,n}^-<
r_{N-1,n}^+<r_{N-2,n}^+<... <r_{p_0+1,n}^+,
$$
and
\[
\lb{gev1}
\g_{0,n}\ss\g_{1,n}\ss...\ss\g_{p_0-1,n}\ss\g_{p_0,n},\ \ \ 
\g_{N-1,n}\ss...\ss\g_{p_0+2,n}\ss\g_{p_0+1,n},
\]
\[ 
\lb{gev2}
 G_{n}=\g_{N-1,n}\cap \g_{0,n}.
\]

Due to \er{gev2} the even gaps are associated with
$f_0^+\approx 1$ and $f_{N-1}^+\approx 1$ and in these two cases
\er{gev2} gives
\[
f_0^+=\cos a-\sin^2 a<1,\qqq f_{N-1}^+=\cos (a-{\pi\/N})-\sin^2(a-{\pi\/N})<1.
\]
Thus the even gaps $G_{n}(a)=(E_{n}^-(a),E_{n}^+(a))$ are determined by
$$
F(E_{n}^\pm(a))=\cos a_+-\sin^2 a_+<1,\qqq a_+=\max \{a,{\pi\/N}-a\}
\in (0,{\pi\/N}).
$$
Thus Lemma \ref{TLas} yields $G_{n}(a)\ne \es$ and \er{T7-4}, in particular,  $|G_n|\to \iy $ as $n\to \iy$.

iii) We consider the  case of  {\bf odd} gaps $\g_{k,n}=(r_{k,n}^-,r_{k,n}^+)$, i.e., $n\in \Z$ is odd. For 
$\l\in \g_{k,n}$ we have 
\[
 F(\l)+s_k^2<-|c_k|\qqq  \lra  \qqq
F(\l)<f _k^-=-|c_k|-s_k^2\in(-{5\/4},-1),\qq k\in \ol{N-1}.
\]
From Theorem \ref{T6} we deduce that these resonances (i.e.,the zeros of $R_k$)  are zeros of Eq. $F(\l)=f_k^-\in [-{5\/4},-1),\l\in \C$. Moreover, by Theorem \ref{T4}, they are real and simple.
We get 
\[
\lb{r1}
\l_{0,n}^- < r_{1,n}^-<r_{2,n}^-<...<r_{p_1,n}^-<r_{p_1,n}^+<r_{p_1-1,n}^+<...
<r_{1,n}^+<\l_{0,n}^+,
\]
\[
\lb{r2}
r_{p_0,n}^-<r_{p_0-1,n}^-<...<r_{p_1+2,n}^-<r_{p_1+1,n}^-
<r_{p_1+1,n}^+<r_{p_1+2,n}^+<...<r_{p_0-1,n}^+<r_{p_0,n}^+,
\]
and
\[
\lb{r3}
r_{p_0+1,n}^-<r_{p_0+2,n}^-<...<r_{p_2,n}^-<r_{p_2,n}^+<r_{p_2-1,n}^+<...
<r_{p_0+1,n}^+,
\]
\[
\lb{r4}
r_{N-1,n}^-<r_{N-2,n}^-<...<r_{p_2+1,n}^-<r_{p_2+1,n}^+<r_{p_2+2,n}^+<...
<r_{N-2,n}^+<r_{N-1,n}^+,
\]
where ${j\pi\/3}-a-{\pi p_j\/N}\in (0,1), j=1,2$ for some integer $p_j\ge 0$.
Thus we obtain
\[
\lb{god1}
\g_{p_1,n}\ss\g_{p_1-1,n}\ss...\ss\g_{1,n}\ss\g_{0,n},\qqq
\g_{p_1+1,n}\ss\g_{p_1+2,n}\ss...\ss\g_{p_0-1,n}\ss\g_{p_0,n},
\]
\[
\lb{god2}
\g_{p_2,n}\ss\g_{p_2+1,n}\ss...\ss\g_{p_0+2,n}\ss\g_{p_0,n},\qqq 
\g_{p_2+1,n}\ss\g_{p_2+2,n}\ss...\ss\g_{N-2,n}\ss\g_{N-1,n},
\]
\[ 
\lb{god3}
 G_{n}=\g_{p_1,n}\cap \g_{p_1+1,n}\cap\g_{p_2,n}\cap \g_{p_2+1,n}.
\]
Then we obtain the odd gaps by 
$G_{n}(a)=(E_{n}^-(a),E_{n}^+(a))$, where $E_{n}^\pm(a)$ are determined by
\[
\lb{logare}
F(E_{n}^\pm(a))=\min \{f_{p_1}^-,f_{p_1+1}^-,f_{p_2}^-,f_{p_2+1}^-\},
\qqq 
f_k^-=-|c_k|-s_k^2.
\]

We will specify \er{logare}.

A) If  $p_1={N\/3}-1\in \Z$, then we get $p_2=2p_1+1$ and
$$
a+{\pi p_1\/N}={\pi\/3}-\ve_1,\qqq \ve_1={\pi\/N}-a\in (0,{\pi\/N})
\qqq a+{\pi p_2\/N}={2\pi\/3}-\ve_1,
$$ 
and 
\[
c_{p_1} =\cos ({\pi\/3}-\ve_1)\ne {1\/2},\qqq 
c_{p_1+1} =\cos ({\pi\/3}-\ve_1+{\pi \/N})\ne {1\/2},
\]
\[
c_{p_2} =\cos ({2\pi\/3}-\ve_1)\ne -{1\/2},\qqq
c_{p_2+1} =\cos ({2\pi\/3}-\ve_1+{\pi \/N})\ne -{1\/2}.
\]
Then due to \er{logare} and Lemma \ref{TLas}  all odd gaps are open $|G_n|>0$ and $|G_n|\to \iy $ as $n\to \iy$.

B) If  $p_2={2N\/3}-1\in \Z, {N\/3}\notin \Z$, then $p_2$ is even and we get $p_1=p_2/2={N\/3}-{1\/2}$ and
$$
a+{\pi p_2\/N}={2\pi\/3}-\ve_1,\qqq
a+{\pi p_1\/N}={\pi\/3}+(a-{\pi\/2N}) 
$$
which give
\[
c_{p_1} =\cos ({\pi\/3}+(a-{\pi\/2N}))\approx {1\/2},\qqq 
c_{p_1+1} =\cos ({\pi\/3}+(a+{\pi\/2N}))\ne {1\/2},
\]
\[
c_{p_2} =\cos ({2\pi\/3}-\ve_1)\ne -{1\/2},\qqq
c_{p_2+1} =\cos ({2\pi\/3}-\ve_1+{\pi \/N})\ne -{1\/2}.
\]
Then due to \er{logare} and Lemma \ref{TLas} we obtain:

if $a={\pi\/2N}$, then all odd gaps are $|G_n|\ge 0$ and $|G_n|\to 0$ as $n\to \iy$. 

If $a\ne {\pi\/2N}$, then all odd gaps are open $|G_n|>0$ and $|G_n|\to \iy $ as $n\to \iy$.

C) Let ${2N\/3}\notin \Z$. 
Let ${jN\/3}=k_j+\wt k_j,\qq \wt k_j\in (0,1), k_j\in \Z,j=1,2$. We have
$$
a+{\pi k_j\/N}={j\pi\/3}+w_j,\qqq w_j\ev a-{\pi \wt k_j\/N}\in (-{\pi \/N},{\pi \/N}),
$$
which give
\[
c_{k_j}=\cos ({j\pi\/3}+w_j),\qqq c_{k_j\pm 1}=\cos ({j\pi\/3}+w_j\pm {\pi \/N})\ne -(-1)^j{1\/2}.
\]

If $a\in \wt A=\{{\pi \wt k_1\/N},{\pi \wt k_2\/N}\}$, then $c_{k_1}={1\/2}, $ or $c_{k_2}=-{1\/2}$. Thus
 we obtain the odd gaps $G_{n}(a)=(E_{n}^-(a),E_{n}^+(a))$, where $E_{m}^\pm(a)$ are determined by
$$
F(E_{m}^\pm(a))=\min \{f_{p_1},f_{p_1-1},f_{p_2},f_{p_2-1}\}=-{5\/4},
\qqq 
$$
and by Lemma \ref{TLas}, $|G_n|\ge 0$ and  $|G_n|\to 0$ as $n\to \iy$.

If $a\notin \wt A$, then we obtain the odd gaps $G_{n}(a)=(E_{n}^-(a),E_{n}^+(a))$, where $E_{m}^\pm(a)$ are determined by \er{logare}:
$$
F(E_{m}^\pm(a))=\min \{f_{p_1},f_{p_1-1},f_{p_2},f_{p_2-1}\}>-{5\/4},
\qqq 
$$
Thus by Lemma \ref{TLas}, $|G_n|>0$ and  $|G_n|\to \iy$ as $n\to \iy$.

iv) In our case $c_k\ne 0$ for all $k=1,..,N$. Then each 
spectral band $\s_{k,n}\ne \es$ for $H_k(a)$  and due to 
\er{r1}-\er{r4}, \er{god1}-\er{god3} we have \er{T7-5}.

v) Let  $a_m={\pi\/2}-{\pi m\/N}\in [0,{\pi\/N}]$ for some $k=m\in \ol N, N\ge 1$. Then $c_m=0$ and in this case by Theorem \ref{T3},
instead of $\s_{m,n}$ we have a flat band $\wt\s_n\{\wt\l_n\}$.
All other spectral bands $\s_{k,n}\ne \es$ for $k\ne m$
and due to 
\er{r1}-\er{r4}, \er{god1}-\er{god3} we have \er{T7-6}.

Using  Theorem \ref{T6} iii), we obtain \er{T7-7}. 
\BBox

\section{Appendix: the direct integral }
\setcounter{equation}{0}

We shortly recall the well known results about 
the properties of point spectrum and the absence of singular continuous spectrum (see \cite{GN}). 

\begin{theorem}\lb{TA} 
For each $(k,a)\in \ol N\ts \R$ the following identities hold:
\[
\lb{TA-1}
\s(H_k)=\s_{ac}(H_k)\cup \s_\iy(H_k),
\]
where the set $\s_\iy(H_k)$ is discrete and  does not have accumulation points,
\[
\lb{TA-2}
\s(H_k)= \s_\iy(H_k),\qqq if\qq c_k=0,
\]
\[
\lb{TA-3}
\s_{ac}(H_k)=\{\l\in \R: {F(\l)+s_k^2\/c_k}\in [-1,1] \},\qqq  if\qq 
c_k\ne 0.
\]
\end{theorem}
\no {\bf Proof.} Recall that for $N=1$  the fundamental subgraph $\G_0$
is given by 
$$
\G_0=\cup_{j=0}^2\G_{0,j,1},\qq
\G_\o=\{r=r_\o^0+te_\o,\  t\in [0,1]\}, \qq r_\o^0,r_\o^1=r_\o^0+e_\o\in \R^3,\qq |e_\o|=1,
$$
where $\G_{0,0}=\G_{0,0,1}$ is  a''vertical`` edge;
$\G_{0,1}=\G_{0,1,1}$ and $\G_{0,2}=\G_{0,2,1}$ are edges with  positive and negative projections on the vector $(0,0,1)\in \R^3$, see Fig. \ref{fig1}.
We have $\G^{(1)}=\cup_{n\in \Z}\G_n$, where $\G_n=\cup_{j=0}^2\G_{n,j}$.
A function $f(x),x\in \G^{(1)}$  has the form $f(x)=f_n(s)$ for $x, s\in \G_n$.
We identify $f_n$ on $\G_n$ with a function on $\G_0$
by using the local coordinate $x=r_\o+te_\o, t\in [0,1]$
and \er{grac}.


Define the space $\gH=\int_{[0,2\pi)}\os \gH_p{dp\/2\pi}$,
where $\gH_p=L^2(\G_0)$. Introduce the unitary operator $U: L^2(\G^{(1)})\to \gH$ and the operator $A_k$  by
$$
(U f)(p)=\sum_{n\in Z}e^{inp}f_n,\quad f_n=(f_{n,j})_{(n,j)\in\Z\ts\Z_n}\in L^2(\G_0), \qq p\in[0,2\pi),
$$ 
$$
A_k=UH_k U^{-1}=\int_{[0,2\pi)}\!\!\!\!\os \ A_k(p){dp\/2\pi},
$$
where an operator $A_k(p)$ on the graph $\G_0$ acts
in the Hilbert space $L^2(\G_0)=\sum_{j\in \Z_3} \os L^2(\G_{0j})$. 
Acting on the edge $\G_{0j}$, $A_k(p)$ is the ordinary differential operator given by 
\[
(A_k(p)f)_j(t)=-f_j''(t)+q(t)f_j(t),\qqq t\in [0,1],
\]
where $f_j, f_j'' \in L^2(\G_{0j}),\ \ j=0,1,2;   q\in L^2(0,1)$ and 
$f\in \gD(A_k(p)$ satisfies 

{\bf The  Kirchhoff  Boundary Conditions}
\[
\lb{aK0}
f_{0}(1)=f_{1}(0)=e^{i a}s^k f_{2}(1),\qqq
e^{ip}f_{0}(0)=e^{i a}f_{1}(1)=f_{2}(0), 
\]
\[
\lb{aK1}
-f'_{0}(1)+f'_{1}(0)-e^{i a}s^k f'_{2}(1)=0,\qqq e^{ip}f_{0}'(0)-e^{ia}f'_{1}(1)+f'_{2}(0)=0.
\]
Remark that \er{aK0}, \er{aK1} follow from \er{1K0}, \er{1K1}, where we used: $f_j=f_{0,j}, j=0,1,2$  and  $f_{n+1,0}(0)=e^{ip}f_{n,0}(0), f_{n+1,0}'(0)=e^{ip}f_{n,0}'(0)$.

Hence $A_k=\int_{[0,2\pi)}\oplus A_k(p)dp$, where
 the operator $A_k(p)$ acts on the finite graph $\G_0$ containing only $3$ edges. Hence $A_k(p)$ has a discrete spectrum and denote
the increasing sequence of the eigenvalues of $A_k(p)$ by 
$E_{k,n}(p), n\ge 1$. If for some $n\ge 1$ the eigenvalue $E_{k,n}(p)=\const$ for all $p\in [0,2\pi]$, then this $E_{k,n}(p)$ is
an eigenvalue of $A_k$ of infinite multiplicity. 
The well-known arguments (see \cite{GN}) give \er{TA-1}.

We solve the equation $A_k(p)f=\l f$.
Recall that any solution $y$ of the equation $-y''+qy=\l y$ satisfies
\[
\lb{a1}
y(t)=w_ty(0)+{\vp_t\/\vp_1}y(1),\qqq w_t=\vt_t-{\vp_t\/\vp_1},\qq t\in[0,1].
\]
Let $x=f_{0}(0),y=f_{0}(1)$. Then using \er{aK0} we obtain
\[
\lb{a22}
f_{0}(t)=w_tx+{\vp_t\/\vp_1}y,\qq
f_1(t)=w_ty+{\vp_t\/\vp_1}e^{i (p-a)}x,
\qq
f_2(t)=e^{ip}w_tx+{\vp_t\/\vp_1}e^{-ia}s^{-k}y.
\]
The substitution   \er{a22} into \er{aK0} gives
\[
-w_1'x-{\vp_1'\/\vp_1}y+w_0'y+{1\/\vp_1}e^{i (p-a)}x-e^{i a}s^k(
e^{ip}w_1'x+{\vp_1'\/\vp_1}e^{-ia}s^{-k}y)=0,
\]
and using $w_0'=-{\vt_1\/\vp_1}, w_1'=-{1\/\vp_1}$ we get
the first Eq.
\[
\lb{ae1}
x(1+e^{i (p-a)}+e^{i (p+a)}s^k)-y(2\D+\vt_1)=0.
\]
We determine the  second one.
The substitution \er{a22}  into \er{aK1} gives
$$
e^{ip}(w_0'x+{1\/\vp_1}y)-e^{ia}(w_1'y+{\vp_1'\/\vp_1}e^{i (p-a)}x)
+(e^{ip}w_0'x+{1\/\vp_1}e^{-ia}s^{-k}y)=0,
$$
which yields the second Eq.
\[
\lb{ae2}
-xe^{ip}(2\D+\vt_1)+y(e^{ip}+e^{ia}+e^{-ia}s^{-k})=0.
\]
The corresponding determinant of the systems \er{ae1}, \er{ae2} has the form
$$
Q=(2\D+\vt_1)(2\D+\vp_1')-(e^{-ip}+e^{-ia}+e^{ia}s^{k})(e^{ip}+e^{ia}+e^{-ia}s^{-k})=0.
$$
Using the identities 
$$
(2\D+\vt_1)(2\D+\vp_1')=8\D^2+\vt_1'\vp_1+1,
$$
$$
(e^{-ip}+e^{-ia}+e^{ia}s^{k})(e^{ip}+e^{ia}+e^{-ia}s^{-k})=
1+4c_k^2+4c_k\cos (p+{\pi k\/N}),
$$
we obtain
\[
Q=F+s_k^2=c_k\cos (p+{\pi k\/N}),\ p\in [0,2\pi).
\]
Thus if $c_k=0$, then $Q=0$, which yields \er{TA-2}.
If $c_k\ne 0$, then we obtain \er{TA-3}.

Note that it is possible to define a modified Lyapunov function
$\wt F_k$ by the identity $\wt F_k(\l,a)={F(\l)+s_k^2\/c_k}=\cos (p+{\pi k\/N}),p\in \R$. This function is entire. But in this case the periodic eigenvalues are defined by the Eq. ${F+s_k^2\/c_k}=\cos {\pi k\/N}$, i.e. at $p=2\pi n$. We assume that we obtain such entire functions only for the zigzag graphs.
\BBox

 \no {\bf Acknowledgments.}
Evgeny Korotyaev was partly supported by DFG project BR691/23-1.
The various parts of this paper were written at the Mittag-Leffler Institute, Stockholm  and in the Erwin Schr\"odinger Institute for Mathematical Physics, Vienna, the first author is grateful to the Institutes for the hospitality.

\end{document}